\documentclass[10pt]{gen-j-l}

\swapnumbers

\usepackage{amssymb,amsfonts}
\usepackage{euscript}

\usepackage[ps,matrix,arrow,curve,cmtip]{xy}

\numberwithin{equation}{section}

\makeatletter
  \let\c@subsection\c@equation
\makeatother

\theoremstyle{plain}   

\newtheorem{thm}[subsection]{Theorem}
\newtheorem{prop}[subsection]{Proposition}
\newtheorem{cor}[subsection]{Corollary}
\newtheorem{lemma}[subsection]{Lemma}

\theoremstyle{remark}
\newtheorem{rem}[subsection]{Remark}    
\newtheorem{exam}[subsection]{Example}

\newtheorem{defn}[subsection]{Definition}

\theoremstyle{plain}

\DeclareMathOperator{\id}{id}
\DeclareMathOperator{\colim}{colim}

\DeclareMathOperator{\Ker}{Ker}

\newcommand{\op}{{\operatorname{op}}}

\newcommand{\Aut}{{\operatorname{Aut}}}
\newcommand{\End}{{\operatorname{End}}}

\newcommand{\alg}[1]{{#1\text{-}\mathrm{alg}}}

\newcommand{\ra}{\rightarrow}

\newcommand{\xra}{\xrightarrow}

\newcommand{\sSet}{{\operatorname{\EuScript{S}}}}
\newcommand{\Top}{{\operatorname{\EuScript{T}}}}

\DeclareMathOperator{\ho}{Ho}
\DeclareMathOperator{\hocolim}{hocolim}

\newcommand{\len}[1]{\lvert#1\rvert}

\newcommand{\powser}[1]{[\![#1]\!]}

\newcommand{\F}{\mathbb{F}}
\newcommand{\Z}{\mathbb{Z}}
\newcommand{\N}{\mathbb{N}}
\newcommand{\R}{\mathbb{R}}
\newcommand{\Q}{\mathbb{Q}}

\newcommand{\point}{{\operatorname{pt}}}
\DeclareMathOperator{\map}{map}

\newcommand{\sm}{\wedge} 
\newcommand{\CP}{\mathbb{CP}}

\newcommand{\red}[1]{\widetilde{#1}}
\newcommand{\pt}{\operatorname{pt}}

\newcommand{\dfn}{\textbf}

\def\noloc{\;{:}\,}

\def\defeq{\overset{\mathrm{def}}=}

\newcommand{\forcepar}{\mbox{}\par}

\title[Logarithmic Cohomology Operation]{The units of a ring spectrum and a
  logarithmic cohomology operation}
\author{Charles Rezk}
\date{ \today}
\address{Department of Mathematics \\
University of Illinois at Urbana-Champaign \\ 
Urbana IL, 61820}
\email{rezk@math.uiuc.edu}
\thanks{The author was partially supported by the National Science
  Foundation.} 

\newcommand{\W}{\mathbb{W}}

\newcommand{\m}{\mathfrak{m}}
\newcommand{\level}{\operatorname{level}}
\renewcommand{\O}{\mathcal{O}}
\newcommand{\I}{\mathcal{I}}

\newcommand{\gbinom}[2]{\genfrac{[}{]}{0pt}{}{#1}{#2}}
\newcommand{\gl}{\operatorname{gl}}

\newcommand{\BK}{\Phi}

\newcommand{\Si}{\Sigma^\infty}
\newcommand{\Sip}{\Sigma^\infty_+}
\newcommand{\Oi}{\Omega^\infty}
\newcommand{\Oim}{\Omega^\infty} 

\newcommand{\cts}{\mathrm{cts}}

\newcommand{\fib}{\mathcal{F}}
\renewcommand{\sSet}{\mathrm{sSet}}
\renewcommand{\Top}{\mathrm{Top}}
\newcommand{\Sing}{\operatorname{Sing}}

\newcommand{\Spectra}{\mathrm{Spectra}}
\newcommand{\Spaces}{\mathrm{Spaces}}
\newcommand{\BFSpectra}{\mathrm{BFSpectra}}

\newcommand{\Lder}{\mathbf{L}}
\newcommand{\Rder}{\mathbf{R}}

\newcommand{\tP}{\widetilde{P}}
\newcommand{\hP}{\widehat{P}}

\newcommand{\oper}[1]{\mathrm{op}_{#1}}

\newcommand{\grplike}{\mathrm{grp}}
\newcommand{\hecke}{\mathcal{H}}

\newcommand{\fset}[1]{\underline{#1}}
\newcommand{\ch}[1]{#1_0^{\sm}}

\newcommand{\bk}{\varphi}
\newcommand{\tbk}{\widetilde{\varphi}}
\newcommand{\tQ}{\widetilde{Q}}

\begin{document}

\begin{abstract}
We construct a ``logarithmic'' cohomology operation on Morava
$E$-theory, which is a homomorphism defined on the multiplicative
group of invertible elements in the ring $E^0(K)$ of a space $K$.  We
obtain a formula for this map in terms of the action of Hecke operators
on Morava $E$-theory.  Our formula is closely related to that for an
Euler factor of the Hecke $L$-function of an automorphic form.
\end{abstract}

\maketitle



\section{Introduction}

Recall that if $R$ is a commutative ring, then the set $R^\times
\subset R$ of invertible elements of $R$ is naturally an abelian group
under multiplication.  This construction is a functor from commutative
rings to abelian groups.  
In general, there is no obvious relation between the additive group of
a ring $R$ and the multiplicative group of units $R^\times$.  However,
under certain circumstances one can define a homomorphism from (a
subgroup of) $R^\times$ to a suitable completion of $R$, e.g., the
natural logarithm $\Q^\times_{>0} \ra \R$, or the $p$-adic
logarithm $(1+p\Z_p)^\times \ra \Z_p$.

The ``logarithmic cohomology operation'' is a homotopy theoretic
analogue of the above, where $R$ is  a commutative $S$-algebra
and ``completion'' is Bousfield localization 
with respect to a Morava $K$-theory.  The purpose of this paper is to
give a formula for the logarithmic operation (in certain contexts) in
terms of power 
operations.
Before giving our results we briefly explain some of the
concepts involved.  

\subsection{Commutative $S$-algebra}

A \dfn{spectrum} is a topological object which represents a
generalized homology and cohomology theory.   A \dfn{commutative
  $S$-algebra} is a spectrum equipped with a commutative
multiplication; such a spectrum gives rise to a cohomology theory
which has a commutative product, as well as power operations.  
Just as any
ordinary commutative ring is an algebra over the ring $\Z$ of
integers, so a commutative $S$-algebra is an algebra over the
\dfn{sphere spectrum} $S$.

The definition of commutative $S$-algebra is rather technical; it is
the result of more than twenty years of effort, by many people.  
There are in fact 
several different models of commutative ring spectra; the commutative
$S$-algebra in 
the sense of \cite{ekmm}, or a symmetric commutative ring spectrum in
the sense of \cite{hovey-shipley-smith-symmetric-spectra}, or some
other equivalent model.  These models are equivalent, in the sense
that they have equivalent homotopy theories; for the purpose of
stating results it does not matter which model we use.

\subsection{Power operations}
A spectrum $R$ which admits the structure of a commutative ring up to
homotopy gives rise to a cohomology theory
$X\mapsto R^*(X)$  taking values in graded commutative
rings.  The structure of commutative $S$-algebra on $R$ is much stronger than
this; it provides not just a ring structure on homotopy groups, but
also \dfn{power operations}, which encode ``higher 
commutativity''.
Let $\Sigma_m$ denote the symmetric group on $m$
letters, and let $B\Sigma_m$ denote its classifying space.
If $R$ is a commutative $S$-algebra, there are 
natural maps
$$P_m\colon R^0(X) \ra R^0(B\Sigma_m \times X),$$
with the property that the composite of $P_m$ with restriction
along an inclusion $\{*\}\times X\ra B\Sigma_m \times X$ is the $m$th
power map $\alpha\mapsto \alpha^m$ on $R^0(X)$.   In other words, for a
commutative $S$-algebra, the $m$th power map is just one of a family
of maps parameterized (in some sense) by the space
$B\Sigma_m$.  An exposition of power operations and their properties
can be found in \cite{bmms-h-infinity-ring-spectra}.  

For suitable $R$, one can construct natural homomorphisms of the
form
$$\rho\colon R^0(B\Sigma_m \times X) \ra D\otimes_{R^0} R^0(X),$$
where $D$ is an $R^0$-algebra, and so get a cohomology operation of
the form 
$$\op\colon R^0(X) \xra{P_m} R^0(B\Sigma_m\times X) \xra{\rho}
D\otimes_{R^0} R^0(X).$$
Such functions $\op$ are  what are usually called power operations.

For instance, the Eilenberg-MacLane spectrum $HR$ associated to
an ordinary commutative ring $R$ is the spectrum which represents
ordinary cohomology;  when $R=\F_p$, power operations are
the Steenrod operations.  Topological $K$-theory spectra, both real
and complex, admit a number of power operations, including the
exterior power operations 
$\lambda^k$ and the Adams operations $\psi^k$.  Other theories of
interest include some elliptic cohomology 
theories, including the spectrum of topological modular forms
\cite{hopkins-icm-2002}; bordism theories, including the spectrum $MU$
of complex bordism.

\subsection{Formal groups and isogenies}
Recall that a multiplicative cohomology theory is \dfn{complex
orientable} if $R^*(\CP^\infty)$ is the ring of
functions on a one-dimensional commutative formal group.  (In this
paper, all formal groups are commutative and one-dimensional.)

If $R$ is a commutative $S$-algebra whose associated cohomology theory
is complex orientable with formal group $G$, and $\op$ is a power
operation on $R$ which is 
a ring homomorphism as above, then 
$$R^0(\CP^\infty) \ra D\otimes_{R^0}R^0(\CP^\infty)$$
is an homomorphism $i^*G\ra G$ of formal groups; here
$i\colon R^0\ra D$ is a map of rings, and $i^*G$ is the formal group
obtained by extension of scalars along $i$.  For example, complex
$K$-theory is complex orientable, and 
$K^0(\CP^\infty)$ is the ring of functions on the formal
multiplicative group $\hat{G}_m$;  the Adams 
operation $\psi^k$  corresponds to the $k$-th power
map $\hat{G}_m\ra \hat{G}_m$.

The philosophy is that power operations (of degree $m$) on a complex
orientable commutative $S$-algebra $R$ should be parameterized by a suitable
family of isogenies (of degree $m$) to the associated formal group.
This philosophy is best understood in the case of Morava $E$-theories,
which we now turn to.  

\subsection{Power operations on Morava $E$-theory}

Fix $1\leq n<\infty$ and a prime $p$.  Let $k$ be a perfect field of
characteristic $p$, and $\Gamma_0$ a height $n$ formal group over
$k$.  Such a formal group admits a Lubin-Tate universal deformation
\cite{lubin-tate-univ-def}, which  
is a formal group $\Gamma$ defined over a ring 
$$\O\approx \W k\powser{u_1,\dots,u_{n-1}};$$
$Wk$ is the ring of $p$-typical Witt vectors on $k$.
\dfn{Morava $E$-theory} is a $2$-periodic complex orientable
cohomology theory with $\Gamma$ as its associated formal group; thus
$\pi_*E\approx \O[u,u^{-1}]$ with $\O$ 
in degree $0$ and $u$ in degree $2$.  
The Hopkins-Miller theorem (see \cite{goerss-hopkins-moduli-spaces},
\cite{richter-robinson-gamma-homology-group-algebras}) states that
Morava $E$-theories admit a canonical structure of 
commutative $S$-algebra.

Power operations for Morava $E$-theories were constructed by Ando
\cite{ando-power-operations}; see also
\cite{ando-hopkins-strickland-h-infinity}.  These operations are
parameterized by  
level structures on the associated Lubin-Tate universal deformation.
To each finite subgroup $A$ of the
infinite torsion group $\Lambda^*\approx (\Q_p/\Z_p)^n$ is associated
a natural ring homomorphism
$\psi_A\colon E^0X\ra D\otimes_{E^0}E^0X$, where $D$ is the
$E_0$-algebra representing a full Drinfel'd 
level structure $f\colon \Lambda^*\ra i^*\Gamma$; the ring $D$ was
introduced 
into homotopy 
theory in
\cite{hopkins-kuhn-ravenel}.  The associated isogeny $i^*G\ra G$ has
as its kernel the subgroup generated by the divisor of the image of
$f|_A$ in $i^*G$.  

An expression in
the $\psi_A$'s which is invariant under the action of the automorphism
group of $\Lambda^*$ descends to a map $E^0X\ra E^0X$ (see
\S\ref{subsec-interpretation-by-hecke} below). 
When $n=1$ and $E$ is
$p$-adic $K$-theory, then $\psi_A=\psi^{p^r}$ for $A\approx
\Z/p^r\subseteq \Lambda^*\approx \Q_p/\Z_p$.  A precise definition of
the operations $\psi_A$ is given in \S\ref{subsec-ando-power-operations}.

\subsection{Units of a commutative ring spectrum}

To a
commutative $S$-algebra $R$ is associated a spectrum $\gl_1(R)$,
which is analogous to the units of a commutative ring (see
\S\ref{sec-units}).  The  
$0$-space of the spectrum 
$\gl_1(R)$ is denoted $GL_1(R)$, and it is equivalent up to
weak equivalence with subspace of $\Oi R$. 

Write $H^q(X;E)\defeq E^q(X)$ for a space $X$ and a spectrum
$E$.  Then $\gl_1(R)$ gives a generalized cohomology theory which, in
degree $q=0$ is given by
$$H^0(X;\gl_1(R)) \approx (R^0(X))^\times.$$
The higher homotopy groups for the spectrum $\gl_1(R)$ are given by
$$\pi_q\gl_1(R)=
\widetilde{H}^0(S^q;\gl_1(R)) \approx (1+\widetilde{R}^0(S^q))^\times
\subset (R^0(S^q))^\times,\qquad q>0.$$
In particular, there is an isomorphism of groups $\pi_q\gl_1(R)\approx
\pi_qR$ for $q>0$, defined by ``$1+x\mapsto x$''.  This isomorphism of
homotopy groups is induced by the inclusion $GL_1(R)\ra \Oi R$ of
spaces, but \emph{not} in general by a map of spectra.

The main interest in the cohomology theory based on $\gl_1(R)$ is that
the group $H^1(X;\gl_1(R))$ contains the obstruction to the
$R$-orientability of vector bundles over $X$, according to the theory
of \cite{may-e-infinity-book}.  The present paper was motivated by one
particular application: the construction of a $MO\langle
8\rangle$-orientation for the topological modular forms spectrum.
This application will appear in  joint work with Matt Ando and
Mike Hopkins.

\subsection{$K(n)$ localization}

Let $F$ be a homology theory.  \dfn{Bousfield $F$-localization}
consists of a functor $L_F$ on the homotopy category of spectra, and a
natural map $\iota_X\colon X\ra L_F X$ for each spectrum $X$, such
that $\iota_X$ is the initial example of a map of spectra out of $X$
which is an $F_*$-homology isomorphism
\cite{bousfield-localization-spectra}.  Distinct homology theories 
may give rise to isomorphic Bousfield localizations, in which case
they are called \dfn{Bousfield equivalent}.

Given a Morava $E$-theory spectrum, there is a an associated ``residue
field'' $F$,  a spectrum formed by killing the sequence of
generators of the ideal $\mathfrak{m}=(p,u_1,\dots,u_{n-1})$ in
$\pi_0E$, so that $\pi_*F\approx 
k[u,u^{-1}]$.  The spectrum $F$ is not a commutative
$S$-algebra, although it is a ring spectrum up to homotopy.  

The Bousfield class of $F$ depends only on the prime
$p$ and the height $n$ of the formal group of $E$, and this is the
same as the Bousfield class of the closely related Morava $K$-theory
spectrum $K(n)$.  (The spectrum $F$ is isomorphic to a finite direct sum of
suspensions of $K(n)$.)  As is standard, we will write $L_{K(n)}$ for
the localization functor associated to any of these Bousfield
equivalent theories. 

In many respects, $K(n)$-localization behaves like completion with 
respect to the ideal $\mathfrak{m}\subset \O$.  In particular,
$K(n)$-localization allows us to define a modification of the homology
functor $E_*$ associated to a Morava-$E$ 
theory, called 
the \dfn{completed $E$-homology} $E^{\sm}_*$ and defined by 
$$E^{\sm}_*(X) \defeq \pi_*(L_{K(n)}(X\sm E)).$$
This functor takes values in complete $E_*$-modules; if $E_*X$ is
a finitely generated $E_*$-module, then $E^{\sm}_*(X)\approx
(E_*(X))^{\sm}_{\mathfrak{m}}$.  See
\cite[\S8]{hovey-strickland-morava-k-theories} for a discussion of
$K(n)$-localization and completed homology.

\subsection{The logarithmic cohomology operation}

For each commutative $S$-algebra $R$, there is a
natural family of ``logarithm'' maps from $\gl_1(R)$ to various
``completions'' of $R$.   
For each prime $p$ and $n\geq 1$, there
exists a natural 
map
$$\ell_{n,p}\colon \gl_1(R) \ra L_{K(n)}R.$$
This map is defined using the construction due to 
Bousfield and Kuhn \cite{bousfield-uniqueness-of-infinite-deloopings},
\cite{kuhn-morava-k-infinity-loop}, which is a functor $\Phi_n$
from spaces to 
spectra, with the property that $\Phi_n \Oi (X)\approx L_{K(n)}X$ for
any spectrum $X$.  If $R$ is a commutative $S$-algebra, the spaces
$\Oi\gl_1(R)$ and $\Oi R$ have weakly equivalent 
basepoint components, and so the Bousfield-Kuhn construction gives an
equivalence $L_{K(n)}\gl_1(R) \approx L_{K(n)}R$ of spectra.  The map
$\ell_{n,p}$ is the composite 
$$\gl_1(R)\xra{\iota_{\gl_1(R)}} L_{K(n)} \gl_1(R) \approx L_{K(n)}
R.$$
The construction of $\ell_{n,p}$ is described in detail  in
\S\ref{sec-construction-logarithm}. 

The map $\ell_{n,p}$ gives a natural transformation of cohomology
theories, and thus for any space $X$ a group
homomorphism 
$$\ell_{n,p}\colon (R^0X)^\times \ra (L_{K(n)}R)^0(X),$$  
natural as $X$ varies over spaces and $R$
varies over commutative $S$-algebras.

The purpose of this paper is the computation of this
``logarithmic'' map in terms of power operations, when $R$ is 
a reasonable $K(n)$-local commutative $S$-algebra.   It is convenient
to consider the cases $n=1$ and $n>1$ separately, though the proof for
$n=1$ is really a corollary of the general case. 

\subsection{The logarithm for $K(1)$-local spectra}
\label{subsec-logarithm-for-k1-local}

Let $p$ be a prime, and let $R$ be a $K(1)$-local commutative
$S$-algebra, satisfying the following technical condition: the kernel of
$\pi_0 L_{K(1)}S \ra \pi_0 R$ contains the torsion subgroup of $\pi_0
L_{K(1)}S$.  This condition is always satisfied if $p>2$ (since $\pi_0
L_{K(1)}S$ is torsion free for odd $p$), and is
satisfied at
all primes when $R$ is the $p$-completion of the periodic complex or real
$K$-theory spectra, or if $R$ is the $K(1)$-localization of the
spectrum of topological modular forms.

Such a ring $R$ admits canonical cohomology operations $\psi$ and
$\theta$, such that (in particular) $\psi$ is a ring homomorphism, and
for $x\in R^0X$,  
$$\psi(x)=x^p+p\theta(x).$$
(See \S\ref{sec-k1-local}.)
When $R$ is the $p$-completion of real or complex
$K$-theory,  then $\psi$ is the classical $p$th Adams operation
$\psi^p$.   

\begin{thm}\label{thm-k1-local-case}
Let $p$ be any prime, and let $R$ be a $K(1)$-local commutative
$S$-algebra, satisfying the technical condition above.  For a finite
complex 
$X$, the logarithm  
$\ell_{1,p}\colon (R^0X)^\times \ra R^0X$ is given by the 
infinite series
$$\ell_{1,p}(x) = \sum_{k=1}^\infty (-1)^k
\frac{p^{k-1}}{k} \left(\frac{\theta(x)}{x^p}\right)^k,$$
which converges $p$-adically for any invertible $x$, and so is a
well-defined expression.  
\end{thm}
Note that the series can be formally
rewritten as
$$\ell_{1,p}(x)=\frac{1}{p}\log\left(
\frac{1}{1+p\,\theta(x)/x^p}\right) = 
\frac{1}{p}\log\left(\frac{x^p}{\psi(x)}\right),$$
and that this new expression is still meaningful, up to $p$-torsion.
Since $x^p/\psi(x)\equiv 1\mod p$ for invertible $x$, this can be written
$$\ell_{1,p}(x)=(\id-\frac{1}{p}\psi)(\log(x)),$$
which is meaningful (up to $p$-torsion) when $x-1$ is nilpotent.  If
$x=1+\epsilon$ with $\epsilon^2=0$, then the formula of
\eqref{thm-k1-local-case} becomes
$$\ell_{1,p}(1+\epsilon) = \epsilon-\theta(\epsilon) =
\epsilon-\frac{1}{p}\psi(\epsilon).$$ 

The proof of \eqref{thm-k1-local-case} is given in
\S\ref{sec-k1-local}, as a corollary of
\eqref{thm-log-formula-lubin-tate-case} below. 

\subsection{The logarithm for Morava $E$-theory}

For general $n\geq1$, we give a result for Morava $E$-theory, in
terms of power operations.
We have
\begin{thm}\label{thm-log-formula-lubin-tate-case}
Let $p$ be any prime, $n\geq 1$, and let $E$ be a Morava $E$-theory
associated to a height $n$ formal group law over a perfect
field of characteristic $p$.
Then the logarithm $\ell_{n,p}\colon (E^0X)^\times \ra E^0X$ is
given by
$$\ell_{n,p}(x)= \sum_{k=1}^\infty (-1)^{k-1}\frac{p^{k-1}}{k} M(x)^k
= \frac{1}{p}\log\left(1+p\cdot M(x)\right),$$
where $M\colon E^0X\ra E^0X$ is the unique cohomology operation such that
$$1+p\cdot M(x) =\prod_{j=0}^n \biggl(\prod_{\substack{
A\subset \Lambda^*[p] \\ \len{A}=p^j }} \psi_A(x)
\biggr)^{(-1)^jp^{(j-1)(j-2)/2}}.$$ 
Here $\Lambda^*[p]\subseteq \Lambda^*$ denotes the kernel of
multiplication by $p$ on $\Lambda^*$.
\end{thm}

In the case when $n=1$, then $1+p\cdot M(x)= x^p/\psi_{\Z/p}(x)$, and
thus we recover the 
result for $K$-theory.  

\subsection{Interpretation in terms of Hecke operators}
\label{subsec-interpretation-by-hecke}

Define formal expressions $T_{j,p}$ for $0\leq j\leq n$ and $T(p^k)$ for
$k\geq 0$ by
$$T_{j,p}\defeq \frac{1}{p^j}\sum_{\substack{A\subseteq \Lambda^*[p]
\\ \len{A}=p^j}}\psi_A \qquad\text{and}\qquad
T(p^k)=\frac{1}{p^k}\sum_{\substack{A\subseteq \Lambda^* \\
\len{A}=p^k}} \psi_A;$$
these give rise to well defined additive operations $E^0X\ra
p^{-1}E^0X$, which we call \dfn{Hecke operators}.  The $T(p^k)$
operators can be 
expressed as polynomials in the $T_{j,p}$ operators, and there is
an action of the commutative ring
$\hecke_{n,p}=\Z[T_{1,p},\dots,T_{n,p}]$ on $p^{-1}E^0X$.  (An account
of this action is sketched in \S\ref{sec-hecke-operators}.)

Formally, we can rewrite the expression of
\eqref{thm-log-formula-lubin-tate-case} as
$$\ell_{n,p}(x) = F_{1}(\log x)$$
(using the fact that $\log$ takes products to sums, and that the
$\psi_A$'s are ring homomorphisms), where
$$F_{X} = \sum_{j=0}^n (-1)^j p^{j(j-1)/2} T_{j,p}\,\cdot X^j \in
\hecke_{n,p}[X].$$ 
In particular, if $x=1+\epsilon \in (E^0X)^\times$ with
$\epsilon^2=0$, then 
$\ell_{n,p}(1+\epsilon)=F_{1}(\epsilon)$ (up to torsion).

The formal operator inverse of $F_X$ is
$$F_{X}^{-1} = \sum_{k=0}^\infty T(p^k)\,\cdot X^k \in
\hecke_{n,p}\powser{X},$$ 
and both these expressions appear in the theory of automorphic forms.
Namely, if $f$ is an eigenvector of the action of the algebra
$\hecke_{n,p}$ on a space of automorphic forms, then
$F_{p^{-s}}^{-1}f=L_p(s;f)f$, where 
$L_p(s;f)=\sum a_k p^{-ks}$ is the $p$-th Euler factor of the Hecke
$L$-function of $f$, and $F_{p^{-s}}f=(1/L_p(s;f))f$, so that
$1/L_p(s;f)$ is a polynomial of degree $n$ in $p^{-s}$.
(See \cite[3.21]{shimura-introduction-automorphic-forms}.)  

It is notable that this expression from the theory of $L$-functions
arises naturally from a purely topological construction; it came as a
suprise to the author, and he still has no good explanation for it.
It is also
significant for the application to elliptic cohomology; in the
presence of an elliptic curve, these Hecke operators coincide with the
classical action of Hecke operators on modular forms.  

\subsection{Structure of the proof}

The proof of \eqref{thm-log-formula-lubin-tate-case} falls naturally
into two parts.  

\textit{First part.}   The logarithm is equal to a certain power
operation in $E$-cohomology, corresponding to a particular element
$v\in\ch{E}{\Oi S}$, the completed $E$-homology of the $0$-th space of
the sphere spectrum \eqref{thm-log-is-power-operation}.  Furthermore,
the element $v$ (the ``logarithmic element'') is completely
characterized by certain algebraic 
properties \eqref{thm-log-comes-from-logarithmic-elt}.  The proof of
the first part comprises
\S\S\ref{sec-bk-operator-formula-basepoint}--\ref{sec-logarithmic-elt}.  

\textit{Second part.}  An element $v\in
\ch{E}{\Oi S}$ is constructed, and shown to be a logarithmic element
\eqref{prop-construction-of-logarithmic-elt}.  The explicit form of
$v$ gives the formula of \eqref{thm-log-formula-lubin-tate-case}.  The
proof of the second part comprises
\S\S\ref{sec-level-structures}--\ref{sec-construction-logarithmic-elt}. 

The proof of \eqref{thm-log-formula-lubin-tate-case} is completed in
\S\ref{subsec-proof-of-main-theorem}.

\subsection{Conventions on spaces and spectra}  
\label{subsec-conventions-spaces-and-spectra}

We write $\Spaces$ for a category of ``spaces'' (such as topological
spaces or simplicial sets), and
$\Spaces_*$ for based spaces.  We write $\Spectra$ for any suitable
category of spectra.  Most of this paper takes place in suitable
homotopy categories of spaces or spectra; therefore, we will usually
not specify a particular model.  There are exceptions, namely
\S\ref{sec-units} (where we refer to the model of \cite{lms}) and
\S\ref{sec-simplicial-functors} (where we use the model of
\cite{bousfield-friedlander}). 

By \emph{commutative $S$-algebra}, we mean any suitable category of commutitive
ring objects in spectra (e.g., the model of \cite{ekmm}, or any
equivalent model, see
\cite{mandell-may-equivariant-orthogonal-spectra-s-modules}).  In 
\S\ref{sec-units}, we will use the particular model of algebras
over the linear isometries operad, in the sense of \cite{lms}.

There are pairs of adjoint functors with units and counits
$$({-})_+\colon \ho\Spaces \rightleftarrows \ho\Spaces_* \noloc
({-})_-, \quad p_K\colon K\ra (K_+)_-, \quad q_K\colon (K_-)_+ \ra
K,$$
where $K_+\approx K\amalg \point$, and $K_-\approx K$ with the
basepoint forgotten,
$$\Si \colon \ho\Spaces_* \rightleftarrows \ho\Spectra\noloc \Oi, \quad
\eta_K\colon K\ra \Oi\Si K, \quad \epsilon_X\colon \Si\Oi X\ra X,$$
and
$$\Sip\colon \ho\Spaces \rightleftarrows \ho\Spectra\noloc \Oim, \quad
\eta^+_K\colon K\ra \Oi\Sip K, \quad \epsilon^+_X \colon \Sip\Oi X\ra
X,$$ 
so that $\Sip K\approx \Si(K_+)$ and $\Oim X\approx (\Oi X)_-$.
Also, note that 
$\Si(X\sm Y) \approx \Si X\sm \Si Y$, while $\Sip(X\times Y) \approx
\Sip X\sm \Sip Y$; we will use this identifications often.

\subsection{Conventions on
localization}\label{subsec-conventions-on-localization}

When $n\geq1$ and the prime $p$ are fixed, we write  $L=L_{K(n)}$ for
the Bousfield localization functor, and write $\iota_X\colon X\ra LX$
for its coaugmentation.

We make the following convention for the sake of legibility: in
general we do not specify the augmentation 
$\iota_X$.  Thus, if $f\colon X\ra Y$ is a map of spectra, the
notation ``$Lf\colon X\ra LY$'' is understood to denote the composite
of $Lf\colon LX\ra LY$ with the coaugmentation $\iota_X\colon X\ra
LX$.  Note that little information is lost, since $Lf\colon LX\ra LY$ is
in fact the unique factorization of $X\ra LY$ through $\iota_X$ up to
homotopy. 

Likewise, if $f\colon X\sm Y\ra Z$ is a map, the notation $Lf \colon
X\sm LY \ra LZ$ denotes the unique extension of $X\sm Y\ra LZ$ along
the map $X\sm Y\ra X\sm LY$ (which is a $K(n)$-homology equivalence).

If $R$ is a $K(n)$-local spectrum, we write
$$R_q^{\sm}X\defeq \pi_q L( X\sm R),\qquad R^qX \defeq [X,\Sigma^qR].$$
Both functors $R^q$ and $R_q^{\sm}$ take $K(n)$-homology
isomorphisms to isomorphisms.

\subsection{Acknowledgements}

This work began as a joint project with Paul Goerss and Mike Hopkins.
In particular, the first proof of \eqref{thm-k1-local-case} when $R$
is $p$-completed $K$-theory was proved jointly with them.  The original
proof of this was somewhat different than the one offered here; it
involved an explicit analysis of the Bousfield-Kuhn functor in the 
$K(1)$-local case.  

I would like to acknowledge both Paul and Mike for their assistance at
various points in this project.  I would also like to thank Matt Ando
for many tutorials on power operations and level structures.  I 
would also like to thank Nick Kuhn and Nora Ganter for various
comments which  improved the paper.

\section{The units of a commutative ring spectrum}
\label{sec-units}

In this section, we describe the \emph{units spectrum} of a
structured commutative ring spectrum.    The notion of the units of a
commutative 
$S$-algebra has a long history, paralleling the long history of
constructions of structured ring spectra.  The notion seems to have
arisen from work of Segal (as in \cite{segal-multiplicative-group})
and Waldhausen.  Our 
discussion  of the units spectrum is  based on the construction of
\cite{may-e-infinity-book}, as corrected in
\cite{may-multiplicative-infinite-loop}.   There is
another approach for constructing the units spectrum due to Woolfson
\cite{woolfson-hyper-gamma-spaces}, based on Segal's theory of
$\Gamma$-spaces. 

\subsection{Definition of the units spectrum}

Let $R$ be a commutative ring spectrum in the sense of \cite{lms},
i.e., a spectrum defined on a universe, and equipped with an action of
the linear isometries operad $\mathcal{L}$.
Then $R(0)=\Oim R$, the $0$-space of
the spectrum $R$, is itself an algebra over $\mathcal{L}$.
Let $GL_1(R)\subseteq \Oim R$ denote the
subspace of $\Oim R$ defined by the pullback square 
$$\xymatrix{
{GL_1(R)} \ar[r]^j \ar[d]
& {\Oi R} \ar[d]
\\
{(\pi_0 R)^\times} \ar[r]
& {\pi_0 R}
}$$
We write $\theta\colon GL_1(R)\ra \Oim R$ for the inclusion.
Then $GL_1(R)$ is a \emph{grouplike} $E_\infty$-space, and so by
infinite loop space theory is the $0$-space of a $(-1)$-connective
spectrum, which we denote $\gl_1 R$.  Note that the identity element
of $GL_1(R)$ is \emph{not} the usual basepoint of $\Oi R$.

The construction which associates $R\mapsto \gl_1 R$ defines a functor
$\ho \alg{S}\ra \ho \Spectra$.  (It can be lifted to an honest zig-zag
of functors between underlying model categories; we don't
need this here.)

\begin{exam}
Let $R=S$.  Then $GL_1(S)\approx G$, the monoid of stable self-homotopy
equivalences of the sphere.  
\end{exam}

\begin{exam}
Let $R=HA$, the Eilenberg-Mac~Lane spectrum associated to a
commutative ring $A$.  Then $\gl_1 HA\approx HA^\times$, the
Eilenberg-Mac~Lane spectrum on the group of units in $A$.  
\end{exam}

The spectrum $\gl_1R$  defines a cohomology theory on spaces;
it is 
convenient to write $X\mapsto H^q(X;\gl_1R)$ for the group represented
by homotopy classes of stable maps from $\Sip X$ to $\Sigma^q\gl_1R$.
In general, there 
seems to be no convenient description of these groups in terms of the
cohomology theory $R$, except when $q=0$.
\begin{prop}\label{prop-units-as-cohomology-group}
There is a natural isomorphism of groups
$$H^0(X;\gl_1R) \approx (R^0X)^\times.$$
Furthermore, if $X$ is a pointed and connected space, then this
isomorphism identifies
$$\widetilde{H}^0(X;\gl_1R) \approx (1+\widetilde{R}^0X)^\times.$$
\end{prop}
In particular, if we take $X=S^k$ for $k\geq1$, we obtain isomorphisms
of groups
$$\pi_k(\gl_1R) \approx \widetilde{H}^0(S^k;\gl_1R) \approx
(1+\widetilde{R}^0S^k)^\times \approx \widetilde{R}^0(S^k) \approx
\pi_k R;$$
this uses the isomorphism $(1+\widetilde{R}^0S^k)^\times\approx
\widetilde{R}^0S^k$ defined by $1+\epsilon\mapsto \epsilon$, and is
realized by a map of spaces $GL_1(R)\ra \Oi R$.  

A main motivation for studying the units is their role in the
obstruction to orientations.  For instance, if $V\ra X$ is a spherical
fibration, the obstruction to the existence of an orientation class in
the $R$-cohomology of the Thom space is a certain class $w(V)\in
H^1(X;\gl_1(R))$; see \cite[IV,\S3]{may-e-infinity-book}.

\subsection{A rational logarithm}

In \S\ref{sec-construction-logarithm} we will construct a
``logarithm'' in the $K(n)$-local setting, 
for $n\geq1$.  This construction does not extend to the case of $n=0$,
where ``$K(0)$-local'' means ``rational''.  For completeness, notice
an \textit{ad hoc} logarithm in the rational 
setting; it will not be used elsewhere in the paper.

Let $R$ be a commutative $S$-algebra, and let $R_\Q$ denote the
rationalization of $R$; its homotopy groups are $\pi_n R_\Q\approx
(\pi_n R)\otimes \Q$.  Let $(\gl_1 R)_1$ denote the $0$-connected
cover of the spectrum $\gl_1 R$.  

The group 
$H^0(X;(\gl_1 R)_1)$ is equal to the subgroup of $(R^0X)^\times$ consisting
of classes $\alpha$ which restrict to $1\in R^0(\{x\})$ for each point
$x\in X$.

\begin{prop}
There exists a map $\ell_0\colon (\gl_1 R)_1\ra R_\Q$ of spectra,
unique up to homotopy,
which when evaluated at a space $X$ is a map
$$H^0(X;(\gl_1R)_1)\ra H^0(X;R_\Q)$$
given by the formula
$$\alpha\mapsto \log(\alpha)=\sum_{k=1}^\infty
(-1)^{k-1}\frac{(\alpha-1)^k}{k}.$$ 
\end{prop}
\begin{proof}
The indicated formula is in fact well-defined; convergence
of the series follows because $H^0(X,R_\Q)\approx \lim
H^0(X^{(k)},R_\Q)$, where $X^{(k)}$ is the $k$-skeleton of a
$CW$-approximation to $X$, and because $\alpha-1\in H^0(X,R)$
is nilpotent when restricted to any finite dimensional complex.

Therefore, the indicated formula gives rise to a natural
transformation of functors to abelian groups, and therefore is
represented by a map $\Oi(\gl_1R)_1\ra \Oi R_\Q$ of $H$-spaces.  Since
$R_\Q$ is a rational spectrum, it is straightforward to show that this
map deloops to a map of spectra, unique up to homotopy.
\end{proof}

When $X=S^k$, $k\geq1$, and $\alpha \in
\widetilde{R}^0S^k$, this gives $\ell_0(1+\alpha)=\alpha$; that is,
$\ell_0$ is the ``identity'' on homotopy groups in
dimensions $\geq1$.


\section{The Bousfield-Kuhn functor and the construction of the
logarithm}\label{sec-construction-logarithm}

Fix a prime $p$ and an integer $n\geq1$.  Write $L=L_{K(n)}$ for
localization of spectra with respect to the $n$th Morava $K$-theory,
as in \S\ref{subsec-conventions-on-localization}.

\subsection{The Bousfield-Kuhn functor}

\begin{prop}[Bousfield
  \cite{bousfield-uniqueness-of-infinite-deloopings}, Kuhn 
\cite{kuhn-morava-k-infinity-loop}, Bousfield
\cite{bousfield-telescopic-localization}]   
There exists a functor $\Phi=\Phi_n\colon \Spaces_* \ra \Spectra$ and
a natural weak equivalence of functors $\tau\colon \Phi\circ \Oi
\xra{\sim} L$.  Furthermore, $\Phi(f)$ is a weak equivalence whenever
$f\colon X\ra Y$ induces an isomorphism on $\pi_n$ for all
sufficiently large $n$.  
\end{prop}

That is, $L\colon \Spectra\ra\Spectra$ factors through
$\Oi\colon \Spectra\ra \Spaces_*$, up to homotopy.

\begin{rem}\label{rem-telescopic-localization}
In fact, a stronger result applies.  There is a functor $\Phi^f_n$ and
a natural equivalence $L_{K(n)}^f=\Phi_n^f\circ \Oi$, where
$L_{K(n)}^f$ is Bousfield localization with respect to a
$v_n$-telescope of a type $n$ finite complex.  In fact, the functor
constructed in \cite{bousfield-telescopic-localization} is $\Phi_n^f$,
in which case $\Phi_n=L_{K(n)}\Phi_n^F$.  

In fact, everywhere in this paper where $L_{K(n)}$ and $\Phi_n$
appear, they may be replaced by $L_{K(n)}^f$ and $\Phi_n^f$, including
in the key results \eqref{thm-log-is-power-operation} and
\eqref{thm-log-comes-from-logarithmic-elt}.
\end{rem}

\subsection{The basepoint shift}
\label{subsec-basepoint-shift}

Let $(K,k_0)$ be a pointed space and $X$ spectrum, 
and let $f\colon K\ra \Oim X$ be an \textit{un}based map.  Write
$j(f)\colon K\ra \Oi X$ for the \textit{based} map defined by
$$k \mapsto (jf)(k)= \mu(f(k), \iota(f(k_0))),$$ where $\mu\colon \Oi
X\times \Oi X\ra \Oi X$ and $\iota \colon \Oi X\ra \Oi X$ are the
addition and inverse maps associated to the infinite loop space
structure.  Colloquially, $(jf)(k) = f(k)-f(k_0)$.  In terms of the
cohomology theory represented by $X$, this induces the evident projection
$X^0(K)\ra \widetilde{X}^0(K)\subset X^0K$ to reduced cohomology summand.

\subsection{Construction of the logarithm}

For a commutative $S$-algebra $R$, we write $\theta\colon GL_1(R)\ra
\Oim R$ for  the standard inclusion, as in \S\ref{sec-units}.  It is not a
basepoint preserving map; 
however, $j\theta$ is.

\begin{defn}
Define
$\ell=\ell_{n,p}\colon \gl_1(R) \ra LR$ to be the composite
$$\gl_1(R) \ra L\gl_1(R) \xra[\sim]{\tau} \Phi(GL_1(R))
\xra{\Phi(j\theta)} 
\Phi(\Oi R) \approx LR.$$
\end{defn}
Note that the map $\Phi(j\theta)$ is an equivalence, since $j\theta$ is a
weak equivalence on basepoint components.  

This construction (for fixed $n$ and $p$) gives rise for each space
$K$ and each commutative $S$-algebra $R$ a map
$$\ell\colon (R^0K)^\times \ra (LR)^0K,$$
which is natural in the variables $K$ and $R$.  It factors through the
composite 
$$(R^0K)^\times \ra ((LR)^0K)^\times \xra{\ell} (LLR)^0K\approx
(LR)^0K;$$
thus, when attempting to calculate the effect of $\ell$ on $R$, it
will suffice to assume that $R$ is already $L$-local.

\section{A formula for the Bousfield-Kuhn idempotent operator}
\label{sec-bk-operator-formula-basepoint}

The Bousfield-Kuhn functor produces  an idempotent  operator which
turns unstable maps between infinite loop spaces into infinite loop
maps; the logarithm of \S\ref{sec-construction-logarithm} is the
result of applying the Bousfield-Kuhn idempotent to the inclusion
$\theta\colon GL_1(R)\ra \Oi R$.  To derive a formula for the logarithm,
we will first give a formula for the Bousfield-Kuhn idempotent.  
In this section, we do this for a version of the Bousfield-Kuhn 
operator $\tbk$ which acts on  basepoint preserving maps.  In the next
section, we extend 
this to an operator $\bk$ which acts on arbitrary maps; that form
will apply to the logarithm.

For based spaces $K, L$, we write $[K,L]_*$ for the set of basepoint
preserving maps up to homotopy.   For spectra $X$, $Y$, we write
$\{X,Y\}$ for the the set of maps in the stable homotopy category.

In what follows, we assume that $X$ and $Y$ are spectra, and that $Y$
is an $L$-local spectrum, so that
$\iota_Y\colon Y\ra LY$ is a weak equivalence.

\subsection{The Bousfield-Kuhn operator}

We define an operator $\tbk\colon [\Oi X,\Oi Y]_*
\ra [\Oi X,\Oi Y]_*$ which sends $f\colon \Oi X\ra \Oi Y$ to the map
obtained by applying $\Oi$ to the composite
$$X\xra{\iota_X} LX \approx \BK\Oi X\xra{\BK f} \BK\Oi Y\approx
LY\approx Y.$$  

The operator $\tbk$ has the following properties, where $f,f'\in [\Oi
  X,\Oi Y]_*$:
\begin{enumerate}
\item [(a)]
$\tbk$ is a natural with respect to maps of spectra $g\colon
X\ra X'$, and to maps of $L$-local spectra $h\colon Y\ra Y'$, in the
sense that
$$\tbk(\Oi h\circ f\circ \Oi g) = \Oi h\circ \tbk f\circ \Oi g.$$
\item [(b)] $\tbk$ is 
additive: $\tbk(f+f')=\tbk f+\tbk f'$, where addition is defined using
the infinite loop structure of $\Oi Y$.
\item [(c)]
If $f=\Oi g$, then $\tbk f=f$.
\item [(d)]
$\tbk f$ is an infinite loop map.
\end{enumerate}
In particular, $\tbk^2=\tbk$, and thus 
the group of 
stable maps $\{X,Y\}$ can be identified with a summand of the group of
unstable maps $[\Oi X,\Oi Y]_*$. 

Our approach to  this operation relies
on two facts.  First, the fact that all unstable maps $f\colon \Oi X\ra
\Oi Y$ between infinite loop spaces factor as $f=\Oi\Si f\circ
\eta_{\Oi X}$, where $\Oi\Si f$ is an infinite loop map and
$\eta_{\Oi X}$ is
the unit of the $\Si$--$\Oi$ adjunction.  Therefore, we really only
need to understand the effect of the Bousfield-Kuhn functor on the
``universal example'' of an unstable map out of an infinite loop
space, which is $\eta_{\Oi X}$.  Second, 
if $\Oi X\approx \Oi \Si K$ where $K$ is a based \emph{space}, then
we can translate the problem of understanding the effect of the
Bousfield-Kuhn functor on $\eta_{\Oi\Si K}$ to
that of understanding its effect on $\eta_{\Oi S}$; this is a
consequence of the fact that the Bousfield-Kuhn functor can be
modelled as a simpliciial functor.

\subsection{The natural transformation $\lambda$}

Let $\lambda_X\colon X\ra L\Si\Oi X$ be the map defined by
$$X\xra{\iota_X} LX \approx \BK\Oi X \xra {\BK \eta_{\Oi X}}
\BK\Oi\Si\Oi X \approx L\Si\Oi X;$$
$\lambda$ is a natural transformation between functors on
the homotopy category of spectra.  It is the same natural
transformation considered in
\cite{kuhn-localization-of-andre-quillen-towers}. 

\begin{prop}\label{prop-phi-in-terms-of-lambda}
Let $f\colon \Oi X\ra \Oi Y$ be a based map, and let
$\widetilde{f}\colon \Si \Oi X\ra Y$ denote its stable adjoint.
There is a commutative diagram
$$\xymatrix{
{X} \ar[r]^-{\lambda_X} \ar[d]_{\BK f\circ \iota_X}
& {L\Si\Oi X} \ar[dl]^{L\widetilde{f}}
\\
{LY\approx Y}
}$$
so that $\tbk f = \Oi(L\widetilde{f}\circ \lambda_X)$.
\end{prop}
\begin{proof}
Observe that $f$ is equal to the
composite
$$\Oi X\xra{\eta_{\Oi X}} \Oi\Si \Oi X \xra{\Oi \widetilde{f}} \Oi
Y.$$ 
Apply $\BK$ to this diagram.
\end{proof}
\begin{cor}
The transformation $\lambda$ is a section of $L\epsilon$.  That is,
$L\epsilon_X\circ \lambda_X=\iota_X\colon X\ra LX$.
\end{cor}
\begin{proof}
Set $f=\Oi \iota_X$ in \eqref{prop-phi-in-terms-of-lambda}.
\end{proof}

\subsection{Formula for $\tbk$}

Let $K$ denote an arbitrary based space.  A map $f\in [\Oi X,\Oi Y]_*$
gives rise to a cohomology operation $f_*\colon \red{X}^0(K)\ra \red{Y}^0(K)$ by$f_*(\alpha)= f\circ \alpha$.
We remind the reader that $f_*$ is not necessarily a homomorphism of
abelian groups, although it is the case that $f_*(0)=0$.
Our goal is to calculate the cohomology operation induced by $\tbk f$ 
in terms of that induced by $f$.  The formula we give is in the form
of a composite
$$\red{X}^0(K)\xra{\tP} \red{X}^0(\Oi S\sm K) \xra{f_*} \red{Y}^0(\Oi
S\sm K) \xra{\tQ} \red{Y}^0(K),$$
where $\tP$ and $\tQ$ are certain functions which we define now.

Given $\alpha\in \red{X}^0(K)$, represented by a map $a\colon \Si K\ra
X$, write $\tP\alpha\in \red{X}^0(\Oi S\sm K)$ for the class represented by 
$$\Si(\Oi S\sm K) \approx \Si\Oi S\sm \Si K \xra{\epsilon_S\sm a} S\sm
X\approx X.$$

If $Y$ is an $L$-local spectrum, then we define a natural map
$$\tQ\colon \red{Y}^0(\Oi S\sm K) \ra \red{Y}^0(K)$$
as follows: represent a class $\alpha\in \red{Y}^0 (\Oi S\sm K)$ by a map
$a\colon \Si(\Oi S\sm K)\ra Y$, and let $\tQ(\alpha)$ be the class
represented by 
$$\Si K\approx S\sm \Si K \xra{\lambda_S\sm 1} L\Si\Oi S\sm \Si K \ra
L\Si( \Oi S\sm K) \xra{La} LY\approx Y.$$

The main result of this section is.
\begin{prop}\label{prop-tbk-identity}
We have an identity $(\tbk f)_* = \tQ\circ f_*\circ \tP$ of operations
$\red{X}^0(K)\ra \red{Y}^0(K)$.  
\end{prop}

\subsection{The natural transformation $\delta$}
\label{subsec-nat-transf-delta}

For a based space $K$ 
and a spectrum $Z$, let $\delta=\delta_{Z,K}\colon  \Oi Z\sm K\ra
\Oi(Z\sm \Si K)$ be the 
map adjoint to 
$$\Si(\Oi Z\sm K) \approx \Si \Oi Z\sm \Si K \xra{\epsilon_Z \sm 
1_{\Si K}} Z\sm \Si K.$$
\begin{lemma}\label{lemma-delta-and-p}
Given $\alpha\in \red{X}^0(K)$ represented by a stable map $a\colon
\Si K\ra X$, the element $\tP(\alpha)$ is represented by the stable
map $\Oi a\circ \delta_{S,K}\colon \Oi S\sm K\ra \Oi X$.
\end{lemma}
\begin{proof}
The composite
$$ \Oi S\sm K\xra{\delta} \Oi(S\sm \Si K) \xra{\Oi a} \Oi X$$
is adjoint to 
$$\Si(\Oi S\sm K)\approx \Si\Oi S\sm \Si K \xra{\epsilon_S \sm 1_{\Si
K}} S\sm \Si K \xra{a} X,$$
which represents $\tP(\alpha)$.
\end{proof}

\subsection{The proof of \eqref{prop-tbk-identity}}

\begin{proof}[Proof of \eqref{prop-tbk-identity}]
Let $\alpha\in \red{X}^0(K)$, represented by a spectrum map $a\colon \Si K\ra
X$.  Let $b\colon \Oi S \sm K\ra \Oi X$ be the map of based spaces
which represents $\tP \alpha\in \red{X}^0(\Oi S\sm K)$; it is adjoit
to the stable map $\epsilon_S\sm a\colon \Si\Oi S\sm \Si K\ra S\sm
X\approx X$.  Let
$\widetilde{f}\colon \Si\Oi X\ra Y$ be the map adjoint to $f\colon \Oi
X\ra \Oi Y$.
We will refer
to the following diagram.
$$\xymatrix{
&& {L\Si ( \Oi S\sm K)}  \ar[d]^{L\Si \delta}
\ar@/^6pc/[dd]^{L\Si b}
\\
{S\sm \Si K} \ar[rr]_-{\lambda_{S\sm \Si K}} \ar[d]_{a}
\ar@(u,l)[urr]^{\lambda_S \sm \iota_{\Si K}}
&& {L\Si\Oi(S\sm \Si K)} \ar[d]^{L\Si\Oi a}
\\
{X} \ar[rr]_-{\lambda_X} \ar[d]_{\BK f\circ \iota_X}
&& {L\Si\Oi X} \ar[dll]^{L\widetilde{f}}
\\
{LY\approx Y}
}$$
We are going to prove that this diagram commutes.  Given this, the
proposition is derived as follows.  Note
that the composite $\BK f\circ \iota_X \circ a$ is precisely the map
representing the class $(\tbk f)_*(\alpha)$.  We claim that the
long composite $S \sm \Si K\ra Y$ around the outer edge of the diagram
is a map representing $(\tQ\circ f_*\circ\tP)(\alpha)$.  The  composite
$\widetilde{f}\circ \Si b\colon \Si(\Oi S\sm K) \ra Y$ is adjoint to
$f\circ b$, which represents $f_*(\tP\alpha)$, and so it is clear from the
definition of $\tQ$ that the long composite in fact the desired class.

To show that the diagram commutes, we need to check the commutativity
of each of four subdiagrams.  The central square commutes because
$\lambda$ defines a natural transformation $1\ra L\Si\Oi$.  The bottom
triangle commutes by \eqref{prop-phi-in-terms-of-lambda}.

The commutativity 
of the right-hand triangle of the diagram follows from
\eqref{lemma-delta-and-p}. 

We defer the proof of the commutativity of the upper-left
triangle\eqref{prop-commutative-lambda-with-rho} to
\S\ref{sec-simplicial-functors}. 
\end{proof}

\section{An unbased Bousfield-Kuhn operator}
\label{sec-bk-operator-formula-nobasepoint}

In this section, we define a Bousfield-Kuhn operator on unbased maps,
using the stable basepoint splitting, and derive a formula for it similar to
\eqref{prop-tbk-identity}; from this we will produce the formuala
\eqref{thm-log-comes-from-logarithmic-elt} for the logarithm.

For unbased spaces $K$ and $L$, we write $[K,L]$ for the set of
unbased homotopy classes of (not necessarily basepoint preserving)
maps.

\subsection{The stable basepoint splitting}
\label{subsec-stable-basepoint-splitting}

For a spectrum $Y$ and a \emph{based} space $K$, we consider functions
$$[K,\Oi Y]_* \xra{i} [K,\Oi Y] \xra{j} [K,\Oi Y]_*.$$
The function $i$ is the evident inclusion, while $j$ is the basepoint
shift operator of \S\ref{subsec-basepoint-shift}.  These
operations give rise to the direct sum decomposition
$$[K,\Oi Y] \approx [K,\Oi Y]_* \oplus \pi_0 Y.$$
In particular, $ji=\id$.
Note this direct sum decomposition arises from the stable
splitting $\Si K_+\approx \Si K\vee S$, which is realized by maps
$$\xymatrix{
{\Sip(\point)} \ar@<1ex>[r]^-{\Sip z_K} 
&
{\Sip K} \ar@<1ex>[l]^-{\Sip \pi_K} \ar@<1ex>[r]^-{\Si q_K} 
&
{\Si K} \ar@<1ex>[l]^-{\gamma_K}
},$$
where $z_K\colon \point\ra K$, $\pi_K\colon K\ra \point$, and
$q_K\colon K_+\ra K$ are the evident maps of spaces, and
$\gamma_K=1_{\Sip K}$ is
the stable map with $\Si q_K\circ \gamma_K=1_{\Si K}$ and $\Sip
\pi_K\circ \gamma_K=0$.   
With this notation, the operators $i$ and $j$ are induced by $q_K$ and
$\gamma_K$ respectively.  
We record the relation between $j$ and $\gamma$ in the following
\begin{lemma}\label{lemma-gamma-formula}
Let $K$ be a based space, and $f\colon K\ra \Oi Y$ an unbased map.
Then the based map $jf$ is adjoint to $\Si K\xra{\gamma_K} \Sip K
\xra{\widetilde{f}} Y$, where $\widetilde{f}$ is adjoint to $f$.  In
particular, $\epsilon_Y^+\circ \gamma_{\Oi Y} = \epsilon_Y$, since
$\epsilon_Y^+$ is adjoint to the identity map of $\Oi Y$, which
preserves the basepoint.

Parameterized Version:  Let $K$ be a based space
and $L$ be an unbased space, and let $f\colon K\times  L\ra \Oi Y$ be
an unbased map, with $\widetilde{f}\colon \Sip K\sm \Sip L\ra Y$ its
adjoint.  Then $\widetilde{f}\circ (\gamma_K\sm 1_{\Sip L})$ is
adjoint to  $f-f\circ(z_K\pi_K\times 1_L)\in [K\times L,\Oi Y]$.

In particular, if $f$ is such that $f\circ (z_K\times 1)\colon L\ra \Oi
Y$ is homotopic to the null map, then $\widetilde{f}\circ (\gamma_K\sm
1_{\Sip L})$ is adjoint to $jf$.
\end{lemma}

\subsection{The natural transformation $\lambda^+$}
\label{subsec-nat-transf-lambdaplus}

Let $\lambda^+_X\colon X\ra L\Sip\Oi X$ be the map defined by
$$X \xra{\lambda_X} L\Si\Oi X \xra{L\gamma_{\Oi X}} L\Sip\Oi X.$$

\subsection{The operator $\bk$}
\label{subsec-operator-bk}

Now suppose that $X$ and $Y$ are spectra, and that $Y$ is $L$-local.
We define an operator $\bk\colon [\Oi X,\Oi Y]\ra [\Oi X,\Oi Y]$ on
the set of unbased maps, by $\bk f\defeq (i\circ \tbk\circ j)f$.
The operator $\bk$ is idempotent, and has as its image the set of
infinite-loop maps; it coincides with $\tbk$ on the summand of
basepoint-preserving maps.  We are going to prove a formula for $\bk$
analogous to the one proved for $\tbk$.

For an \emph{unbased} space $K$ we define natural functions
$$P\colon X^0(K) \ra X^0(\Oi S\times K)\qquad
\text{and}\qquad Q \colon Y^0(\Oi S\times K) \ra
Y^0(K),$$ 
as follows.

Given $\alpha\in X^0(K)$, represented by a map $a\colon \Sip K\ra
X$, write $P\alpha\in X^0(\Oi S\times K)$ for the class represented by  
$$\Sip(\Oi S\times K) \approx \Sip\Oi S\sm \Sip K \xra{\epsilon_S^+\sm
a} S\sm X\approx X.$$

If $Y$ is an $L$-local spectrum, then we define a natural map
$$Q\colon Y^0(\Oi S\times K) \ra Y^0(K)$$
as follows: represent a class $\alpha\in Y^0 (\Oi S\times K)$ by a map
$a\colon \Sip(\Oi S\times K)\ra Y$, and let $\tQ(\alpha)$ be the class
represented by 
$$\Sip K\approx S\sm \Sip K \xra{\lambda_S^+\sm 1} L\Sip\Oi S \sm \Sip
K\ra L\Sip(\Oi S\times K) \xra{La} LY\approx Y.$$

\begin{prop}\label{prop-gamma-property}
We have an identity $(\bk f)_* = Q\circ f_*\circ P$ of operations
$X^0(K)\ra Y^0(K)$.  
\end{prop}
\begin{proof}
Let $K$ be an unbased space, and let $\alpha \in X^0(K) \approx
\red{X}^0(K_+)$.  We will show that 
$Q(f_*(P\alpha)) = \tQ((jf)_*(\tP(\alpha))) = (\tbk j f)_*(\alpha)$;
the result follows when we note that there is an identity of
cohomology operations $(i g)_*=g_*$ when $g\in [\Oi
X,\Oi Y]_*$, so that $(\tbk jf)_*=(i\tbk jf)_*=(\bk f)_*$.

We have a direct sum decomposition
$$X^0(\Oi S\times K) \approx
X^0(K)\oplus \red{X}^0(\Oi S\sm K_+) \approx X^0(K)\oplus X^0(\Oi
S\times K,\point\times K),$$
which is produced by smashing $K_+$ with the maps
$$(\pt)_+ \xra{z} (\Oi S)_+ \xra{q} \Oi S \qquad \text{and} \qquad
\pi\colon (\Oi S)_+ \ra (\pt)_+.$$
The remaining projection of
the splitting comes from smashing $\Sip K$ with the stable map
$\gamma=\gamma_{\Oi S} \colon  
\Si \Oi S\ra \Si (\Oi S)_+$.
We claim that with respect to this splitting, the maps in
$$\red{X}^0(K_+) \xra{P} \red{X}^0(K_+)\oplus \red{X}^0(\Oi S\sm K_+)
\xra{f_*} \red{Y}^0(K_+)\oplus \red{Y}^0(\Oi S\sm K_+) \xra{Q}
\red{Y}^0(K_+)$$
satisfy
$$P(\alpha)=(0,\tP(\alpha)),\quad f_*(0,\beta)=(f_*(0),
(jf)_*(\beta)), \quad Q(\alpha,\beta)=\tQ(\beta).$$ 
The proposition will follow immediately.

The formula for $P$ follows, using a standard adjunction argument,
from the facts that $\epsilon_S^+\circ \gamma=\epsilon_S$, by
\eqref{lemma-gamma-formula}, and $\epsilon_S^+ \circ \Si
z=0$.

The formula for $Q$ follows from the facts that $L\Si q\circ
\lambda_S^+=L\Si q \circ L\gamma \circ \lambda_S = \lambda_S$
and $L\Si \pi \circ \lambda_S^+ = L\Si \pi\circ L\gamma \circ
\lambda_S = 0$.

To prove the formula $f_*$, let $b\colon \Oi S\times K\ra \Oi X$ be the
map representing $\beta \in \red{X}^0(\Oi S\sm K_+)\approx
X^0(\Oi S\times K, \point\times K)$; in particular, $b\circ (z\times
1)\colon \point\times K\ra \Oi X$ is the null map, and $b$ can be
taken to be basepoint preserving. Set $f_*(0,\beta)=(x,y)$.  It
is 
clear that $x$ is represented by $f\circ b\circ (z\times
1)=f\circ 0\colon \point\times K\ra \Oi Y$, and so $x=f_*(0)\in
X^0(K)\approx 
\red{X}^0(K_+)$. 
We have that $y$ is represented by the stable map $\widetilde{f}\circ
\Sip b\circ \gamma\sm 1$, which by the parameterized
version of \eqref{lemma-gamma-formula} is adjoint to $j(f\circ b)$.
Since $b$ preserves basepoints, $j(f\circ b)=(jf)\circ b$, which
represents $(jf)_*(\beta)$.

\end{proof}

\subsection{Application to ring spectra}

We now introduce the additional hypothesis that $Y$ is a commutative
ring spectrum; here we only need that $Y$ be a commutative ring up to
homotopy, not a structured ring spectrum.  We write $\ch{Y}(X)\defeq
\pi_0 L(X\sm  Y)$ for the completed homology of a spectrum  $X$ with
respect to  
$Y$.  If $K$ is a based space, we write $\ch{Y}(K)$ for $\ch{Y}(\Sip
K)$. 

Let $v\in \ch{Y}(\Oi S)$ denote the Hurewicz image of $\lambda_S^+\in
\pi_0 L\Oi S$; i.e., the homology class represented by
$$S\xra{1_S\sm \kappa} S\sm Y \xra{\lambda^+_S\sm 1_Y} L\Sip\Oi S\sm Y
\ra L(\Sip \Oi S\sm Y),$$
where $\kappa\colon S\ra Y$ represents the unit of the ring spectrum.

Since $Y$ is a ring spectrum, there is a slant product operation
$$\alpha\otimes \sigma \mapsto \alpha/ \sigma \colon
Y^0(X_1\sm  X_2)\otimes \ch{Y}(X_2) \ra Y^0(X_1)$$ 
defined by 
$$X_1 \xra{1\sm \sigma} X_1\sm L(X_2\sm Y)
\ra L(X_1\sm X_2\sm Y) \xra{L(\alpha\sm 1)} L(Y\sm Y) \ra
LY\approx Y.$$

\begin{prop}
For an $L$-local ring spectrum $Y$, and $\alpha\in Y^0(K)$, we have
$$Q(\alpha)=\alpha/ v.$$
In particular, for $\alpha\in X^0(K)$, we have  $(\bk f)_*(\alpha) = 
f_*(P\alpha)/ v$.  
\end{prop}
\begin{proof}
The second statement is immediate from the first.  The first statement
is straightforward from the definitions.
\end{proof}

Now suppose that $Y=R$ is an $K(n)$-local commutative $S$-algebra, and
that $X=\gl_1(R)$.  Let $\theta\colon \Oi \gl_1(R)\ra \Oi R$ be the
standard inclusion; it corresponds to the ``cohomology operation''
$R^0(K)^\times \ra R^0(K)$ which is the inclusion of the units into
the ring.  The logarithm is $\ell = (\bk \theta)_*$, and thus we have
proved
\begin{thm}\label{thm-log-is-power-operation}
Let $v\in \ch{R}\Sip \Oim S$ denote the Hurewicz image of $\lambda^+_S \in
\pi_0 L\Sip\Oim S$.  Then $\ell(\alpha)=\theta_*(P\alpha)/v$.
\end{thm}
We will often abuse notation by taking $\theta_*$ to be understood, so that we
write $\ell(\alpha)=(P\alpha)/v$.

\section{Simplicial functors}
\label{sec-simplicial-functors}

For a spectrum $X$ and a space $K$, recall (from
\S\ref{subsec-nat-transf-delta}) that $\delta_{X,K}\colon \Oi X\sm
K_+ \ra \Oi(X\sm \Sip K)$ is the map adjoint to
$\epsilon_X\sm\id\colon \Si\Oi X\sm \Sip
K\ra X\sm \Sip K$, where $\epsilon$ is the counit of the adjunction
$(\Si,\Oi)$. 
The purpose of this section is to prove the following.  
\begin{prop}\label{prop-commutative-lambda-with-rho}
For every spectrum $X$ and space $K$, the diagram
$$\xymatrix{
{LX\sm \Sip K} \ar[rr]^-{\lambda_X\sm \id} \ar[d]
&& {L(\Si\Oi X)\sm \Sip K} \ar[d]^{L(\Si \delta_{X,K})}
\\
{L(X\sm \Sip K)} \ar[rr]_-{\lambda_{X\sm \Sip K}}
&& {L\Si\Oi (X\sm \Sip K)}
}$$
commutes in the homotopy category of spectra; the vertical maps
involve the localization augmentation, as described in
\S\ref{subsec-conventions-on-localization}. 
\end{prop}
The idea of the proof is easy to describe: the functors in
question come from \emph{simplicial} functors on an underlying simplicial
model category of spectra, and the natural transformation $\lambda$
on the homotopy category comes from a natural transformation between
these simplicial functors; the vertical maps in the square are certain
natural transformations associated a simplicial functor.  To carry out
the proof, we will choose explicit models for these functors.

\subsection{Simplicial functors}

Recall that if $C$ is a simplicial model category, then for any
$X,Y\in C$ and $K\in\sSet$ there are objects
$$\map_C(X,Y)\in \sSet, \quad X\otimes K\in C, \quad Y^K\in C,$$
which come with isomorphisms
$$\map_C(X\otimes K,Y)\approx \map_C(X,Y^K)\approx
\map_{\sSet}(K,\map_C(X,Y)).$$ 
We will usually write $T_KX\defeq X\otimes K$ and $M^KY\defeq Y^K$;
these are functors $T_K,M_K\colon C\ra C$.  It is standard that these
give rise to derived functors $\Lder T_K, \Rder M_K\colon \ho C\ra \ho
C$, and that furthermore $\Lder T_K$ and $\Rder M_K$ are adjoint to
each other.  See \cite[IX]{goerss-jardine-simplicial-book} for a
full treatment.

A \dfn{simplicial functor} $F\colon C\ra D$ is a functor which is enriched
over simplicial sets, in the sense that for every pair of objects
there is an induced map $\map_C(X,Y)\ra \map_D(FX,FY)$ which on
$0$-simplicies coincides with $F$.  For any such functor, there is a
natural transformation
$$\rho^F_{X,K} \colon F(X)\otimes K\ra F(X\otimes K),$$
which is adjoint to the map 
$$K\ra \map_C(X\otimes \point, X\otimes K) \ra \map_D(F(X\otimes
\point), F(X\otimes K)).$$

A \dfn{simplicial natural transformation} is a natural transformation
$\alpha\colon F\ra G$ such that the two evident maps $\map_C(X,Y)\ra
\map_D(FX,GY)$ sending $f$ to $\alpha_Y\circ Ff$ and $Gf\circ
\alpha_X$ coincide.  Such transformations give rise to commutative
square
$$\xymatrix{
{FX\otimes K} \ar[r]^{\rho^F_{X,K}} \ar[d]_{\alpha_X\otimes 1_K}
& {F(X\otimes K)} \ar[d]^{\alpha_{X\otimes K}}
\\
{GX\otimes K} \ar[r]_{\rho^G_{X,K}}
& {G(X\otimes K)}
}$$

\subsection{Bousfield-Friedlander spectra}

Let $\BFSpectra$ denote the Bousfield-Friedlander model category of
spectra \cite{bousfield-friedlander},
\cite{hovey-shipley-smith-symmetric-spectra},
\cite[X.4]{goerss-jardine-simplicial-book}. 
This category has as objects 
$X=\{X_n\in \sSet_*,f^X_n\colon S^1\sm X_n\ra X_{n+1}\}_{n\geq0}$, where
where $S^1=\Delta[1]/\partial\Delta[1]$, with morphisms $g\colon
X\ra Y$ being sequences $\{g_n\colon
X_n\ra Y_n\}_{n\geq0}$ commuting with the structure maps $f$.  It
simplicial model category.

We note the following functors; all of them are simplicial functors.
\begin{enumerate}
\item [(1)] 
$\Si\colon \sSet_*\ra \BFSpectra$, defined on objects by
$(\Si K)_n = {(S^1)}^{\sm n} \sm K$, with the evident structure maps.
\item [(2)] 
$\Oi \colon \BFSpectra$, defined on objects by $\Oi
X=X_0$.
\item [(3)] 
$\fib\colon \BFSpectra\ra \BFSpectra$, defined on objects by
$$(\fib X)_n = \colim_m \Omega^m\Sing\lvert X_{m+n}\rvert,$$
where $\Sing\colon \sSet_*\rightleftarrows \Top_*\noloc \lvert
{-}\rvert$ are geometric 
realization and singular complex, and $\Omega K=\map_*(S^1,
K)$.  (The functors $\lvert\;\rvert$, $\Sing$, and $\Omega$ are
simplicial functors, and thus so is $\fib$.)
\item [(4)] 
$\Phi\colon \sSet_*\ra \BFSpectra$, the functor defined in
\cite{bousfield-telescopic-localization}; there it is constructed as a
simplicial functor.  
\item [(5)] 
$T_K\colon \sSet_*\ra \sSet$, where $K$ is an unpointed simplicial
set, defined on objects by $T_K(L)=L\sm K_+$.  
\item [(6)] 
$T_K\colon \BFSpectra\ra \BFSpectra$, where $K$ is an unpointed
simplicial set, defined on objects by $T_K(X)_n = X_n\sm K_+$.
\end{enumerate}

A spectrum  $X\in \BFSpectra$ is \emph{cofibrant} if and only if each
structure map
$f^X_n\colon X_n\ra S^1\sm X_{n+1}$ is an inclusion of simplicial
sets.
For all $K\in \sSet_*$, the spectrum $\Si K$ is cofibrant.

A map $f\colon X\ra Y\in \BFSpectra$ between \emph{cofibrant} objects
is a weak equivalence if and only $\fib f$ is a weak equivalence of
simplicial sets in each degree. 
The functor $\fib$ comes equipped with a natural transformation $1\ra
\fib$, which on degree $n$ of a spectrum $X$ is the map adjoint to the
evident map $\lvert X_n\rvert \ra \colim_m \Oi\lvert X_{m+n}\rvert$.

The functors $\Si$, $\Phi$, and both versions of $T_K$ are
\emph{homotopy functors}: they preserve weak equivalences.   The
functor $\fib$ is a homotopy functor on the full subcategory of
cofibrant objects.
The
functor $\Oi$ is not a homotopy functors; however, the composite
$\Oi\fib$ is a homotopy functor on the full subcategory of cofibrant
spectra. 

There is a simplicial natural transformation $\eta\colon \id\ra
\Oi\fib\Si$ of 
functors $\sSet_*\ra \sSet_*$; it is the evident map
$$X\ra \colim_m \Omega^m  \Sing\lvert S^m\sm X\rvert.$$

Define 
$$\lambda\colon \Phi \Oi\fib =\Phi\id\Oi\fib \ra
\Phi \Oi\fib \Si\Oi\fib$$ 
to be the natural simplicial transformation $\Phi\eta\Oi\fib$. 

By the above remarks, both functors are homotopy functors on cofibrant
spectra.  Both functors are simplicial functors.  Therefore, the
square 
$$\xymatrix{
{T_K\Phi \Oi \fib X} \ar[r]^-{T_K\lambda_X} \ar[d]_{\rho}
& {T_K\Phi\Oi \fib\Si\Oi \fib X} \ar[d]^{\rho'}
\\
{\Phi\Oi \fib T_KX} \ar[r]_-{\lambda_{T_KX}}
& {\Phi\Oi\fib\Si\Oi\fib T_KX}
}$$
commutes for all $X$.  Applied to cofibrant $X$, it gives rise to a
commutative diagram in the homotopy category of spectra.  It is clear
that the horizontal arrows are precisely those of
\eqref{prop-commutative-lambda-with-rho}.  It remains
to identify the arrows labelled $\rho$ and $\rho'$.  They can be
factored 
$$\rho = \rho^{\Phi\Oi \fib}\quad \text{and} \quad \rho' =
(\Phi\Oi\fib) (\Si) \rho^{\Oi\fib}\circ (\Phi\Oi\fib) \rho^{\Si}(\Oi \fib)
\circ \rho^{\Phi\Oi\fib}(\Si)(\Oi \fib).$$

\begin{prop}
\forcepar
\begin{enumerate}
\item [(a)] 
The natural transformation
$\rho^{\Phi\Oi\fib}\colon T_K\Phi\Oi\fib\ra \Phi\Oi\fib T_K$, on the
homotopy category of spectra, induces the map
$$LX\sm \Sip K\ra L(X\sm \Sip K),$$
which is the unique map up to homotopy compatible with the augmentations
$\iota_X\sm \id\colon X\sm \Sip K\ra LX\sm \Sip K$ and $\iota\colon
X\sm \Sip K\ra L(X\sm \Sip K)$. 
\item [(b)]
The natural transformation
$\rho^{\Si}\colon T_K\Si\ra \Si T_K$, on the homotopy category of
spectra, induces the canonical equivalence 
$$\Si X\sm \Sip K\xra{\sim} \Si(X\sm \Sip K).$$
\item [(c)]
The natural transformation
$\rho^{\Oi\fib}\colon T_K\Oi\fib \ra \Oi\fib T_K$, on the homotopy
category of spectra, induces the map
$$\delta_{X,K}\colon \Oi X\sm \Sip K\ra \Oi(X\sm \Sip K)$$
adjoint to $\epsilon_X\sm \id\colon \Si\Oi X\sm \Sip K\ra X\sm \Sip
K$. 
\end{enumerate}
\end{prop}
\begin{proof}
\begin{enumerate}
\item [(a)] The composite $\Phi\Oi \fib$ is a model for the Bousfield
localization functor on spectra.    
\item [(b)] The map $\rho^{\Si}$ is easily seen to be an isomorphism
on objects in $\BFSpectra$.
\item [(c)] The map $\rho^{\Oi\fib}$ is adjoint to the arrow $\alpha$
in the square
$$\xymatrix{
{T_K\Si\Oi\fib} \ar[r]^-{T_K\epsilon\fib}
\ar[d]_{\rho^{\Si}\Oi\fib}_{\sim}
& {T_K\fib} \ar[d]^{\rho^{\fib}}_{\sim}
\\
{\Si T_K\Oi \fib} \ar[r]_-\alpha
& {\fib T_K}
}$$
which is seen to commute by appeal to the definitions.  The vertical
arrow on the left is an isomorphism, and the vertical arrow on the
right is a weak equivalence when applied to cofibrant spectra.  This
gives the desired result.  
\end{enumerate}
\end{proof}

\section{Power operations and the units spectrum}
\label{sec-power-ops-and-units}

Power operations on a representable functor
on the homotopy category of spaces are defined via the action
of some $E_\infty$-operad on the representing space.  We need to
consider two flavors of $E_\infty$-spaces, and their associated power
operations: 
the ``additive'' $E_\infty$-structure on the $0$-space of any
spectrum (power operations here amount to the theory of transfers),
and the ``multiplicative'' $E_\infty$-structure on the $0$-space of
any commutative $S$-algebra.  The purpose of this section is to show
that both flavors 
are united in the $E_\infty$-space $GL_1(R)$,  which can be viewed as 
``multiplicative'' by the inclusion $GL_1(R)\subseteq \Oi R$, or as
``additive'' by $GL_1(R)\approx \Oi \gl_1(R)$.  

There are some standard references for some of this material, notably
\cite{bmms-h-infinity-ring-spectra}, who frame their power
operations in terms of extended power of spectra.  (Their formulation
works naturaly not just for commutative $S$-algebras (a.k.a.,
$E_\infty$-ring spectra), 
but with the weaker notion of an $H_\infty$-ring spectra, and much of
what we say in this section applies in that context, too.)

\subsection{Power operations associated to $E_\infty$-spaces}

Let $Z$ be an $E_\infty$-space, as in \cite{geometry-iterated-loop}.
That is, for a suitable 
$E_\infty$-operad $\mathcal{E}$, there are structure maps
$\mathcal{E}_m \times_{\Sigma_m} Z^m \ra Z$, satisfying certain
axioms.  Each space $\mathcal{E}_m$ of the operad is a free
contractible $\Sigma_m$-space; we will not refer explicitly to the
operad $\mathcal{E}$, but rather to a structure map
$E\Sigma_m\times_{\Sigma_m} Z^m\ra Z$.

Given a map $f\colon K\ra Z$, we
let $\hP_m(f)\colon E\Sigma_m \times_{\Sigma_m} K^m\ra Z$ denote the
composite $E\Sigma_m\times_{\Sigma_m}K^m\ra E\Sigma_m
\times_{\Sigma_m} Z^m \ra Z$. 
This gives rise to functions $\hP_m$ and $P_m$ on homotopy classes of
maps 
$$\xymatrix{
{[K,Z]} \ar[r]^-{\hP_m} \ar[dr]_{P_m}
& {[E\Sigma_m\times_{\Sigma_m}K^m,Z]} \ar[d]^{\Delta^*_m}
\\
& {[B\Sigma_m\times K,Z]}
}$$
Here $\Delta_m\colon B\Sigma_m\times K\approx
E\Sigma_m\times_{\Sigma_m} K\ra E\Sigma_m\times_{\Sigma_m} K^m$
denotes the diagonal map, and $P_m\defeq \Delta^*_m\circ \hP_m$.

The set $[K,Z]$ is a commutative monoid via the
$E_\infty$-structure of $Z$.  With respect to this, $\hP_m$ and $P_m$
are homomorphisms of commutative monoids.  Furthemore, these
operations are natural with respect to maps of spaces in the
$K$-variable, and maps of $E_\infty$-spaces in the $Z$-variable.

\subsection{Power operations associated to a commutative $S$-algebra} 
\label{subsec-power-ops-commutative-s-alg}

Let $R$ be a commutative
$S$-algebra, and $Z=\Oi R$, viewed as an $E_\infty$-space using the
``multiplicative'' structure (which, we emphasize, is \emph{not} the
same as the ``additive'' structure coming from $Z$ being an infinite
loop space).  Then the constructions of the previous section
specialize to natural power maps $\hP_m$ and $P_m$:
$$\xymatrix{
{R^0K} \ar[r]^-{\hP_m} \ar[dr]_-{P_m}
& {R^0(E\Sigma_m\times_{\Sigma_m} K^m)} \ar[d]^{\Delta^*_m}
\\
& {R^0(B\Sigma_m\times K)}
}$$
These operations are \emph{multiplicative}:
$\hP_m(\alpha\beta)=\hP_m(\alpha)\hP_m(\beta)$ and $\hP_m(1)=1$, and
similarly for $P_m$.
Hence if  $\alpha\in R^0K$ is multiplicatively invertible, then
$\hP_m(\alpha)$ and $P_m(\alpha)$
are also invertible, and so these operations restrict to functions
$(R^0K)^\times \ra (R^0E\Sigma_m\times_{\Sigma_m}K^n)^\times \ra (R^0
B\Sigma_m\times K)^\times$.    

These operations are not in general additive.  Instead, we have
\begin{prop}\label{prop-power-op-of-sum} 
We have
$$\hP_m(\alpha+\beta)=\sum_{i+j=m}
\widetilde{T}_{ij}[P_i(\alpha)P_j(\alpha)]$$
and 
$$P_m(\alpha+\beta)=\sum_{i+j=m} T_{ij}[P_i(\alpha)P_j(\beta)],$$
where $\widetilde{T}_{ij}$ and $T_{ij}$ denote the transfers
associated to the covering maps $E\Sigma_m\times_{\Sigma_i\times
  \Sigma_j}K^m \ra E\Sigma_m\times_{\Sigma_m}K^m$ and $B\Sigma_i\times
B\Sigma_j\ra B\Sigma_m$ respectively.
\end{prop}
\begin{proof}
This is \cite[Lemma 2.1]{bmms-h-infinity-ring-spectra}.
\end{proof}

\subsection{Power operations associated to an infinite loop space, and
transfers} 
\label{subsec-power-op-infinite-loop}

Let $Y$ be a spectrum, and let $Z=\Oi Y$, viewed as an
$E_\infty$-space using the ``additive'' structure.  Then we again have
natural power maps, which we denote $\hP^+_m$ and $P^+_m$:
$$\xymatrix{
{Y^0K} \ar[r]^-{\hP^+_m} \ar[dr]_-{P^+_m}
& {Y^0(E\Sigma_m\times_{\Sigma_m} K^m)} \ar[d]^{\Delta^*_m}
\\
& {Y^0(B\Sigma_m\times K)}
}$$
These operations are \emph{additive}:
$\hP_m^+(\alpha+\beta)=\hP_m^+(\alpha)+\hP_m^+(\beta)$ and $\hP_m^+(0)=0$, and
so likewise for $P_m^+$.  

We recall the well-known relation between such ``power operations''
and the theory of transfers.

\begin{prop}\label{cor-adjoint-of-universal-mth-power} 
  The map $\hP_m^+(\eta^+_K)\colon
  E\Sigma_m\times_{\Sigma_m}K^m\ra \Oim\Sip K$ is adjoint to the
  composite
$$\Sip(E\Sigma_m\times_{\Sigma_m}K^m)\xra{\text{transfer}}
  \Sip(E\Sigma_m\times_{\Sigma_m}(\fset{m}\times K^m))
  \xra{\Sip(\text{proj})} \Sip K,$$
where $\fset{m}$ is a fixed set of size $m$ permuted by $\Sigma_m$.
In particular, $\hP_m^+(\eta^+_*)\colon B\Sigma_m\ra \Oim\Sip S$ is
  adjoint to the composite
$$\Sip B\Sigma_m \xra{\text{transfer}} \Sip B\Sigma_{m-1}
\xra{\Sip(\text{proj})} \Sip(*)\approx S,$$
(take $B\Sigma_{-1}=\varnothing$).
\end{prop}

\subsection{The natural transformation $\delta^+$}
\label{subsec-nat-transf-deltaplus}

Let $Y$ be a spectrum and $K$ an unbased space.  Let
$\delta^+_{Y,K}\colon \Oim Y\times K\ra \Oim(Y\sm \Sip
K)$
denote the map adjoint to 
$$\epsilon^+_Y\sm \id \colon \Sip\Oim Y\sm \Sip K\ra Y\sm
\Sip K.$$
This is a variant of the map $\delta$ defined in
\S\ref{subsec-nat-transf-delta}.

\subsection{A total power operation for infinite loop spaces}
\label{subsec-total-power-op-infinite-loop}

Let $K$ be a space and $Y$  a spectrum.  We define
operations $\hP\colon Y^0K \ra Y^0\Oim\Sip K$, and $P\colon
Y^0K\ra Y^0(\Oim\Sip S\times K)$, as follows.  Given
$\alpha\colon K\ra \Oim 
Y$, let $\hP(\alpha)\defeq \Oim\widetilde{\alpha}$, where
$\widetilde{\alpha}\colon \Sip K\ra Y$ is the adjoint to $\alpha$.
Let $P(\alpha) \defeq \hP(\alpha)\circ \delta^+_{S,K}$.

Since $\Oim Y$ is an infinite loop space, it admits
``additive'' power operations of the type described in
\S\ref{subsec-power-op-infinite-loop}. 
The following lemma says that the two kinds of operations coincide,
via the standard maps 
$$\widetilde{i}_m\defeq \hP_m^+(\eta^+_K)\colon
E\Sigma_m\times_{\Sigma_m}K^m \ra \Oim\Sip K 
\quad\text{and}\quad i_m\defeq P_m^+(\eta_*^+)\colon
B\Sigma_m \ra \Oim S.$$  
\begin{lemma}\label{lemma-operations-coincide}
Let $\alpha\colon K\ra \Oim Y$ be a map.
The diagram
$$\xymatrix{
{B\Sigma_m \times K} \ar[r]^-{\Delta} \ar[d]_{i_m\times
\id_K} 
& {E\Sigma_m \times_{\Sigma_m} K^m} \ar[r]^-{\hP_m^+(\alpha)}
\ar[d]_{\widetilde{i}_m} 
& {\Oim Y}
\\
{\Oim S\times K} \ar[r]_{\delta^+_{S,K}} 
& {\Oim\Sip K} \ar[ur]_{\hP(\alpha)}
}$$
commutes up to homotopy. 
\end{lemma}
\begin{proof}
The adjoint pair $(\Sip,\Oim)$  gives a commutative diagram
$$\xymatrix{
{K} \ar[r]^\alpha \ar[d]_{\eta^+_K}
& {\Oim Y} 
\\
{\Oim\Sip K} \ar[ur]_{\Oim\tilde\alpha}
}$$
The two maps $\alpha$ and $\eta^+_K$ from $K$ give rise to maps
$\hP_m^+(\alpha)$ and $\hP_m^+(\eta^+_K)$ out of
$E\Sigma_m\times_{\Sigma_m} 
K^m$, and the resulting triangle which appears in the 
statement of the proposition commutes, by the naturality of $\hP_m^+$
with respect to $E_\infty$-maps.

The commutativity of the left-hand square comes from the fact that the
map $\hP_m^+(\eta_X^+)\colon 
E\Sigma_m\times_{\Sigma_m} X^m \ra \Oim\Sip X$ is a natural transformation of
functors, and can in fact be realized as a \emph{toplological}
natural transformation, by taking $E\Sigma_m$ to be the $m$-th space
of the little  cubes operad.
\end{proof}

The example we are interested in is $Y=\gl_1(R)$.  In this
case, ``classical'' power operations are just the standard power
operations in $R$-cohomology, restricted to units.  

\subsection{The universal example of the ``infinite-loop'' total power
operation}

The ``infinite-loop'' power
operation $\hP\colon [K,\Oim Y]\ra [\Oim \Sip K,\Oim Y]$ of the
previous section is completely
determined by its restriction to a ``canonical'' class, namely the
map $\eta^+_K\colon K\ra \Oim\Sip K$.  This leads to the following.
\begin{prop}\label{prop-alt-formula-inf-loop-power-ops}
Given a map $\alpha\colon K\ra \Oim Y$, we have that 
\begin{enumerate}
\item [(a)] $\hP(\alpha)\colon \Oim\Sip K\ra \Oim Y$ is adjoint to the
composite $\Sip\Oim\Sip K \xra{\epsilon_{\Sip K}^+} \Sip K
\xra{\widetilde{\alpha}} Y$.
\item [(b)] $P(\alpha)\colon \Oim S\times K\ra \Oim Y$ is adjoint to the
composite 
$$\Sip\Oim S\sm \Sip K \xra{\epsilon_{S}^+\sm\id_{\Sip K}} S\sm \Sip K
\approx \Sip K \xra{\widetilde{\alpha}} Y.$$
\end{enumerate}
Here $\widetilde{\alpha}$ denotes the adjoint to $\alpha$.  In
particular, the operator $P$ coincides with the one called $P$ in
\S\ref{subsec-operator-bk}. 
\end{prop}
\begin{proof}
Consider the diagram
$$\xymatrix{
{\Sip \Oim S\sm \Sip K} \ar[rr]^{\Sip\delta^+_{S,K}}
\ar[drr]_{\epsilon_S^+\sm\id_{\Sip K}}
&& {\Sip\Oim \Sip K} \ar[rr]^{\Sip\Oim\widetilde{\alpha}}
\ar[d]^{\epsilon_{\Sip K}^+} 
&& {\Sip\Oim Y} \ar[d]^{\epsilon_Y^+} 
\\
&& {\Sip K} \ar[rr]_{\widetilde{\alpha}}
&& {Y}
}$$
It is clear that the square commutes (up to homotopy), because
$\epsilon^+$ is a natural transformation, while the triangle commutes
because $\delta_{S,K}^+$ is adjoint to $\epsilon_S^+\sm 1_{\Sip K}$.

$\hP(\alpha)=\Oim \widetilde{\alpha}$ is adjoint to
$\epsilon_Y^+\circ \Sip\Oim\widetilde{\alpha}$, which equals
$\widetilde{\alpha}\circ \epsilon^+_{\Sip K}$ by commutativity of the
diagram.  Thus $P(\alpha)= \hP(\alpha)\circ \delta^+_{S,K}$ is
adjoint to $\widetilde{\alpha}\circ \epsilon^+_{\Sip K}\circ
\Sip\delta^+_{S,K}$, which the diagram shows is equal to
$\widetilde{\alpha}\circ \epsilon_S^+\sm\id_{\Sip K}$.  
\end{proof}

\subsection{Power operations associated to homology classes}
\label{subsec-power-op-of-homology-class}

We now assume that $R$ is a $K(n)$-local commutative $S$-algebra.
From now on, we 
will be interested only in two related kinds of power operations: the
operations $\hP_n$ 
and $P_n$ associated to the mulplicative struture on $\Oi R$
(\S\ref{subsec-power-ops-commutative-s-alg}), and the 
total operations $\hP$ and $P$ associated to the infinite loop space
$GL_1(R)$ (\S\ref{subsec-power-op-infinite-loop}).

Given $u\in \ch{R}B\Sigma_m$, 
define
$\oper{u}\colon R^0K\ra R^0K$ by $\oper{u}(x)\defeq  P_{m}(x)/u$,
using the slant product map ${-}/u\colon R^0(B\Sigma_m \times K)\ra
R^0(K)$.  

Similarly, for a class $u\in \ch{R}\Oim S$ we define $\oper{u}\colon
R^0(K)^\times \ra 
R^0(K)$ by $\oper{u}(\alpha)\defeq P(\alpha)/u$, where we
implicitly regard $P(\alpha)$ as an element of $R^0(\Oi S\times K)$ by the
usual inclusion $R^0(\Oi S\times K)^\times\subseteq R^0(\Oi S\times K)$.
Such 
operations defined using $\Oim S$ coincide with those defined using
$B\Sigma_m$, via the map $i_n=\hP_m(\eta_*^+)$ used in
\eqref{lemma-operations-coincide}.  That is, if $u\in R_0B\Sigma_n$,
and $u'={i_m}_*(u)\in R_0\Oim S$, then $\oper{u}=\oper{u'}$.  

With these definitions, \eqref{thm-log-is-power-operation} becomes
$$\ell(\alpha)=\op_v(\alpha).$$

More generally, if
$D$ is a  flat extension of $\ch{R}(\point)$, then we obtain
operations $\oper{u}\colon R^0K\ra D\otimes_{R}R^0K$
parameterized by elements $u\in D\otimes_R \ch{R}B\Sigma_m$
(respectively, $\oper{u}\colon (R^0K)^\times\ra D\otimes_R R^0K$
parameterized by elements $u\in D\otimes_R \ch{R}\Oi S$).

\section{The structure of the spaces parameterizing power operations}
\label{sec-structure-spaces-param-power-ops}

We summarize here some structure which is relevant to power
operations, and which is used in
\S\S\ref{sec-logarithmic-elt}--\ref{sec-construction-logarithmic-elt}.
Most of what we say in this section is well-known; the summary is
provided mainly to fix notation.
For another summary of this sort of structure, in a 
similar context, the reader is directed to
\cite{strickland-turner-rational-morava-e-theory}.

Let $R$ be a $K(n)$-local commutative
$S$-algebra, and $D$  a flat extention of $\ch{R}(\point)$. Set
$h(X)\defeq D\otimes_{\ch{R}(\point)} \ch{R}(X)$; this is a
multiplicative, homological functor.  Set $h^0(X)\defeq
D\otimes_{\ch{R}}(\point) R^0X$. 
We write $h$ for $h(\point)$.   The functor $h$ admits a K\"unneth map
$\times\colon h(X)\otimes_h h(Y)\ra h(X\times Y)$.   
Let the symbol $M$  denote \emph{either} $B\Sigma$ or $\Oi S$.  

We give below a combined list of structure maps (numbered items)
involving the 
group $h(M)$, and a list of
properties they satisfy (lettered items).  
After the list, we will give the
definitions of each of the structure maps, and sketch the proofs of
each of the properties.   We will also show that in every case in which
the structure maps defined for both $M=B\Sigma$ and $M=\Oi S$, they
commute with the standard map $B\Sigma\ra \Oi S$.  
It may be helpful to point
out here that structures (1) through (7) make $h(\Oi S)$ into a Hopf ring.

\subsection{Summary list of structure maps for $B\Sigma$ and $\Oi S$}
\label{subsec-structure-maps}

The letters $u,v$ denote elements of $h(M)$, while $\alpha$
denotes an element of $R^0(K)$, for an arbitrary space $K$.

\begin{enumerate}
\item [(1)] $h(M)$ is a  module over $h=h(\point)$, such that

\begin{enumerate}
\item [(a)] $\oper{u+v}(\alpha)=\oper{u}(\alpha)+\oper{v}(\alpha)$,
and  $\oper{cu}(\alpha)=c\cdot\oper{u}(\alpha),$ for $c\in h(\point)$.
\end{enumerate}

\item [(2)] There is a distinguished class $1\in h(M)$, and

\item [(3)] a product $u\otimes
  v\mapsto u\cdot   
  v\colon   h(M)\otimes_{h} h(M)\ra h(M)$, such that

\begin{enumerate}
\item [(b)] ``$\cdot$'' is  associative and commutative, with
unit $1$, and 

\item [(c)] $\oper{1}(\alpha)=1$ and $\oper{u\cdot
v}(\alpha)=\oper{u}(\alpha)\oper{v}(\alpha)$.
\end{enumerate}
\end{enumerate}
Structures (1), (2), and (3) can be summarized: $h(M)$ is a commutative
$h$-algebra, and for $\alpha\in h^0(K)$, the map
$u\mapsto \oper{u}(\alpha)\colon h^0(M)\ra h^0(K)$ is a map of
$h$-algebras. 

\begin{enumerate}
\item [(4)] There is a distinguished class $[1]\in h(M)$, and

\item [(5)] a product $u\otimes   v\mapsto u\circ   v\colon
h(M)\otimes_{h} h(M)\ra h(M)$, such that

\begin{enumerate}
\item [(d)] ``$\circ$'' is associative and commutative, with unit
$[1]$, and

\item [(e)] $\op_{[1]}=\id$.
\end{enumerate}

\item [(6)] There is an $h$-module map $\pi\colon h(M)\ra h$, such that

\begin{enumerate}
\item [(g)] $\pi(1)=1$ and $\pi(u\cdot v)=\pi(u)\pi(v)$, and
\item [(h)]  $\pi(u)=\oper{u}(1)$.  
\end{enumerate}

\item [(7)] There is an $h$-module map $\Delta^\times\colon h(M)\ra
  h(M\times M)$. 

Say that an element $u\in h(M)$ is \dfn{grouplike} if (i)
$\pi(u)=1$, and (ii) $\Delta^\times(u)=u\times u$.  We write
$h(M)^{\grplike}$ for the set of grouplike elements. 

\begin{enumerate}
\item [(i)] If $u\in h(M)$ is grouplike, then $\oper{u}$ is multiplicative:
$\oper{u}(1)=1$ and
$\oper{u}(\alpha\beta)=\oper{u}(\alpha)\oper{u}(\beta)$. 

\item [(j)] If $u$ is grouplike, then we have the identities
$u\circ 1=1$ and $u\circ(v\cdot w)=(u\circ v)\cdot
(u\circ w)$. 
That is, $x\mapsto u\circ x$ is a ring homomorphism if $u$ is grouplike.

\item [(k)] $1\in h(M)^\grplike$, and $u,v\in h(M)^\grplike$ implies
$u\cdot v\in h(M)^\grplike$ and $u\circ v\in h(M)^\grplike$.  Thus,
the set of grouplike elements  
$h(M)^\grplike$ forms a commutative semi-ring, in which ``addition''
is given by 
$u\cdot v$, ``multiplication'' is given by
$u\circ v$, the ``zero'' element is $1$, and the ``one'' element
is $[1]$.  Furthermore, the grouplike elements of $h(\Oi S)$ admit
``additive'' inverses, so that $h(\Oi S)^\grplike$ is not just a
semi-ring but a ring.
\end{enumerate}

\item [(8)] There is an $h$-module map $\tau\colon h(M)\ra h$, such that

\begin{enumerate}
\item [(l)]  $\tau([1])=1$ and $\tau(u\circ v)=\tau(u)\tau(v)$, and

\item [(m)] $\tau(1)=0$ and $\tau(u\cdot v)=
\tau(u)\pi(v)+\pi(u)\tau(v)$. 
\end{enumerate}
In particular, $\tau$ defines a semi-ring homomorphism
$h(M)^\grplike\ra h$.
\end{enumerate}
The following structures are only defined for $M=B\Sigma$.

\begin{enumerate}
\item [(9)] There is an $h$-module map $\zeta\colon h(B\Sigma)\ra
  h(\point)$, such that

\begin{enumerate}
\item [(n)] $\zeta(1)=1$ and $\zeta(u\cdot v)=\zeta(u)\zeta(v)$, 
\item [(o)] $\zeta([1])=0$ and $\zeta(u\circ
v)=\zeta(u)\pi(v)+\pi(u)\zeta(v)$, and
\item [(p)] $\oper{u}(0)=\zeta(u)$.  
\end{enumerate}

\item [(10)] There is an $h$-module map $\Delta^+\colon h(B\Sigma)\ra
h(B\Sigma\times B\Sigma)$.

Say that $u\in h(B\Sigma)$ is \dfn{primitive} if (i) $\zeta(u)=0$, and (ii)
$\Delta^+(u)=u\times 1+1\times u$.

\begin{enumerate}
\item [(q)] If $u\in h(B\Sigma)$ is primitive,
then $\oper{u}$ is additive:
$\oper{u}(0)=0$ and
$\oper{u}(\alpha+\beta)=\oper{u}(\alpha)+\oper{u}(\beta)$. 
\end{enumerate}
\end{enumerate}

\subsection{Definition of the structure maps}

Recall that $\Sip B\Sigma\approx DS$, the free
commutative $S$-algebra on the $0$-sphere.

\begin{enumerate}
\item [(2)] The class $1\in h(M)$ is the image of the canonical class
under the maps $\point
\approx B\Sigma_0\subset B\Sigma\xra{i} \Oi S$.  Equivalently, it is
induced by the unit map $S\xra{1} DS$ of the ring $DS$.

\item [(3)] The product ``$\cdot$'' is induced on $B\Sigma$ and $\Oi
S$ by maps
$$B\amalg \colon B\Sigma\times B\Sigma\ra B\Sigma\quad \text{and}
\quad H_\amalg\colon \Oi S\times \Oi S\ra \Oi S.$$
Here $\amalg\colon \Sigma\times \Sigma\ra \Sigma$ is the coproduct
functor on finite sets; $\Sip B\amalg$ corresponds to the ring product
$DS\sm DS\ra DS$.  $H_\amalg$ is the ``additive'' $H$-space
product on $\Oi S$, obtained by applying $\Oi$ to the fold map
$S\times S\approx S\vee S\ra S$. 

\item [(4)] The class $[1]\in h(M)$ is the image of the canonical
class under the maps $\point\approx B\Sigma_1\subset B\Sigma \xra{i}
\Oi S$.  Equivalently, it is induced by the ``cannonical'' map $S\ra
DS$. 

\item [(5)] The product ``$\circ$'' induced on $B\Sigma$ and $\Oi S$
by maps 
$$B\times\colon B\Sigma\times B\Sigma\ra B\Sigma \quad\text{and}\quad
H_\times\colon B\Sigma\times B\Sigma\ra B\Sigma.$$
Here $\times \colon \Sigma\times \Sigma\ra \Sigma$ is the cartesian
product functor on finite sets, and $H_\times$ is the
``multiplicative'' $H$-space product on $\Oi S$, and is adjoint to
$\epsilon^+\sm \epsilon^+\colon \Sip\Oi S\sm \Sip\Oi S\ra S\sm
S\approx S$.  

\item [(6)] The map $\pi$ is induced by the
projection map $M\ra \point$.

\item [(7)] The map $\Delta^\times$ is induced by the diagonal map
$M\ra M\times M$.

\item [(8)] The map $\tau$ is induced on $\Sip B\Sigma$ by the composites
$$\Sip B\Sigma_m \xra{\text{transfer}} \Sip B\Sigma_{m-1}
\xra{\Sip(\text{proj})} \Sip(\point)=S,$$
(let $B\Sigma_{-1}=\varnothing$).
On $\Sip\Oi S$ the map $\tau$ is induced by the counit map
$\epsilon^+\colon 
\Sip\Oim S\ra S$.

\item [(9)] The map $\zeta$ is induced on $\Sip B\Sigma$ by maps $j_k\colon
\Sip B\Sigma_n\ra \Sip B\Sigma_0\approx S$, which is the identity if
$k=0$, and null homotopic if $k>0$.  Alternately, $\zeta$ is induced
by the ring map $DS\ra S$ free on the null map $S\ra *\ra S$.

\item [(10)] The map $\Delta^+$ is induced on $\Sip B\Sigma$ by the
transfer maps 
$$\Sip B\Sigma_k \ra \Sip (B\Sigma_i\times B\Sigma_j),\quad i+j=k,$$
associated to the inclusion $\Sigma_i\times \Sigma_j\subset
\Sigma_n$.  Alternately, $\Delta^+$ is induced by the ring map $DS\ra
D(S\vee S)\approx DS\sm DS$ which is free on the pinch map $S\ra S\vee S$.
\end{enumerate}

\subsection{Compatibility of the structure maps}

We need to know that in each of the cases (1) through (8), the two
structure maps are compatable with respect to $B\Sigma\ra \Oi S$.

For (1), (2), (4), (6), and (7), compatibility is clear.

For (3) and (5), compatibility amounts to saying that $B\Sigma \ra \Oi
S$ is a map of ``semi-ring spaces'', which is well-known.

For (8), this is \eqref{cor-adjoint-of-universal-mth-power}.

\subsection{Proof of properties}

\begin{enumerate}
\item [(a)] The slant product $h^0(M\times K)\otimes_h h(M)\ra h^0(K)$
  is $h$-linear.
\item [(b)] $B\amalg$ (resp.~$H_\amalg$) make $M$
  into a commutative and associative $H$-space.
\item [(c)] For $M=B\Sigma$, this follows from the fact that
  $r_{ij}^*\hP_m(\alpha)=\hP_i(\alpha)\times \hP_j(\alpha)$, where
  $i+j=m$ and   $r_{ij}\colon (E\Sigma_i\times_{\Sigma_i} K^i) \times
  (E\Sigma_j\times_{\Sigma_j} K^j)\approx
  E\Sigma_m\times_{\Sigma_i\times \Sigma_j} K^m \ra
  E\Sigma_m\times_{\Sigma_n}K^m$.  Restricting along the diagonal maps
  $\Delta_i$ and $\Delta_j$
  gives the result. 

  For $M=\Oi S$, this follows from the fact that
  $H_\amalg^*\hP(\alpha)=\hP(\alpha)\times \hP(\alpha)$, viewed as a
  class in $h^0(\Oim \Sip K\times \Oim\Sip K)^\times$.  
\item [(d)] $B\times$ (resp.~$H_\times$) make $M$ into a commutative
  and associative $H$-space.
\item [(e)] $P_1\colon h^0(K)\ra h^0(E\Sigma_1\times_{\Sigma_1} K)\approx
  h^0(K)$ is the identity map.
\item [(g)] $M\ra \point$ is a map of $H$-spaces.
\item [(h)] $P\colon h^0(\point)\ra h^0(M\times \point)$ sends $1$ to
  $1$, and $\pi(u)$ equals the slant product of $1$ with $u$.
\item [(i)] If $u$ is grouplike, then the slant product map
  ${-}/u\colon h^0(M\times K)\ra h^0K$ is a map of rings.  Since
  $P\colon h^0(K)\ra h^0(M\times K)$ is mulitplicative, the result
  follows. 
\item [(j)] $M$ is a (semi-)ring space.
\item [(k)] The two products are defined via space-level maps
  $M\times M\ra M$, and so are compatible with diagonal.  
\item [(l)] For $M=B\Sigma$, this is a consequence of the
  ``double-coset formula'' applied to the homotopy pullback square
$$\xymatrix{
{B\Sigma_{k-1}\times B\Sigma_{\ell-1}} \ar[r] \ar[d]
& {B\Sigma_{k\ell-1}} \ar[d]
\\
{B\Sigma_k\times B\Sigma_\ell} \ar[r]^-{B\times}
& {B\Sigma_{k\ell}.}
}$$
For $M=\Oi S$, it follows from the fact that the composite
$$\Sip \Oim S\sm \Sip \Oim S \xra{\Sip H_\times} \Sip\Oim S
\xra{\epsilon^+} S,$$
which induces $u\otimes v\mapsto \tau(u\circ v)$,
equals $\epsilon^+\sm \epsilon^+$.
\item [(m)] For $M=B\Sigma$, this is a consequence of the
  ``double-coset formula'' applied to the homotopy pullback square
$$\xymatrix{
{(B\Sigma_{k-1}\times B\Sigma_\ell) \amalg (B\Sigma_k\times
  B\Sigma_{\ell-1})} \ar[r] \ar[d]
& {B\Sigma_{k+\ell-1}} \ar[d]
\\
{B\Sigma_k \times B\Sigma_\ell} \ar[r]^-{B\amalg}
& {B\Sigma_{k+\ell}.}
}$$
To prove it for $M=\Oi S$, note that for any spectrum $X$ there is a
diagram
$$\xymatrix{
{\Sip\Oim (X\times X)} \ar[rrr]^-{(\Sip\Oim \pi_1, \Sip\Oim \pi_2)}
\ar[drrr]_{\Sip\Oim(\text{fold})}
&&& {\Sip\Oim X\vee \Sip\Oim X} \ar[d]^{\text{fold}}
\ar[r]^-{\epsilon^+\vee \epsilon^+}
& {X\vee X} \ar[d]^{\text{fold}}
\\
&&& {\Sip\Oim X} \ar[r]_{\epsilon^+}
& {X}
}$$
where $\pi_i\colon X\times X\ra X$ denote the projection maps.  This
diagram is commutative up to homotopy.  Evaluating at $X=S$ gives the
desired result.
\item [(n)] This is clear from the observation that $\zeta$ arises
  from a ring map $DS\ra S$.
\item [(o)] $B\times$ maps $B\Sigma_k\times B\Sigma_\ell$ to
  $B\Sigma_0$ only if either $k$ or $\ell$ equals $0$.
\item [(p)] $P\colon h^0(\point)\ra h^0(B\Sigma)\approx \prod_k
  h^0(B\Sigma_k)$ sends $0$ to $(1,0,0,\dots)$.
\item [(q)] This follows from \eqref{prop-power-op-of-sum}.
\end{enumerate}

\section{The logarithmic element}
\label{sec-logarithmic-elt}

In this section, we characterize the element $v\in \ch{R}(\Oi S)$
which appears in the 
statement of
\eqref{thm-log-is-power-operation}.  The main result is
\eqref{thm-log-comes-from-logarithmic-elt}, which states that $v$ is
satisfies certain algebraic 
identities, which characterize it uniquely; i.e., it is the unique
\emph{logarithmic element}.  These algebraic identities in some sense
encode the properties of the idempotent operator $\bk$: namely, that
$\bk$ is the identity on infinite loop maps, and is idempotent.

We note that all the results in this section still hold if we take $L$
to be the telescopic localization functor $L_{K(n)}^f$, rather than
$K(n)$-localization; see \eqref{rem-telescopic-localization}. 

\subsection{Logarithmic elements}

Recall from \S\ref{subsec-structure-maps}(8) the homomorphism
$\tau\colon \ch{R}\Oi S\ra \ch{R}S$
induced by the counit map $\epsilon^+_S\colon \Sip\Oim S \ra
\Sip(\pt)=S$, and from \S\ref{subsec-structure-maps}(5) the product
$\circ \colon \ch{R}\Oi S\otimes_{\ch{R}(\point)} \ch{R}\Oi S\ra
\ch{R}\Oi S$.

A \dfn{logarithmic element} for $R$ is an element $v\in \ch{R}\Oi S$ with
the following two properties: 
\begin{enumerate}
\item [(La)] $\tau(v)=1$ in $\ch{R}S$.
\item [(Lb)] $x\circ v = \tau(x) v$ for all $x\in
\ch{R}\Oi S$.  
\end{enumerate}
\begin{prop}\label{prop-uniqueness-of-logarithmic-element}
There is at most one logarithmic element in $\ch{R}\Oi S$.  A map
$R\ra R'$ of $K(n)$-local commutative
$S$-algebras carries the logarithmic element for $R$ (if
it exists) to the logarithmic element for $R'$.
\end{prop}
\begin{proof}
The second statement is a consequence of the uniqueness of the logarithmic
element.  If both $v$ and $v'$ are logarithmic elements in $\ch{R}\Oi
S$, then
$$v=1\cdot v = \tau(v')\cdot v = v'\circ v = v\circ v' = \tau(v)\cdot
v' = 1\cdot v'=v',$$
using \S\ref{subsec-structure-maps}(b) and (d).
\end{proof}

\begin{thm}\label{thm-log-comes-from-logarithmic-elt}
The element $v$ of \eqref{thm-log-is-power-operation} is a logarithmic
element in $\ch{R}\Oim S$, and hence is the unique logarithmic element.
\end{thm}
In particular, there is a logarithmic element for $R=LS$, the
$K(n)$-localization of the sphere.  By
\eqref{prop-uniqueness-of-logarithmic-element} all logarithmic elements for
$K(n)$-local commutative $S$-algebras are determined by the logarithmic
element for the case $R=LS$.

To give a proof of \eqref{thm-log-comes-from-logarithmic-elt}, we need
some results involving the natural transformations $\lambda^+$
(\S\ref{subsec-nat-transf-lambdaplus}) and $\delta^+$
(\S\ref{subsec-nat-transf-deltaplus}). 

Recall the structure map $\gamma$ of
\S\ref{subsec-stable-basepoint-splitting}.  
\begin{prop}\label{prop-delta-commutes-with-gamma}
The diagram 
$$\xymatrix{
{\Si\Oi Y\sm \Sip K}  \ar[rr]^{\gamma_{\Oi Y}\sm 1_{\Sip K}}
\ar[d]_{\Si \delta_{Y,K}} 
&& {\Sip\Oim Y\sm \Sip K} \ar[d]^{\Sip \delta^+_{Y,K}}
\\
{\Si\Oi(Y\sm \Sip K)} \ar[rr]_{\gamma_{\Oi (Y\sm \Sip K)}}
&& {\Sip\Oim(Y\sm \Sip K)}
}$$
commutes in the homotopy category of spectra. 
\end{prop}
\begin{proof}
This is immediate from
\eqref{lemma-delta-maps-commute-with-projections}, and the stable basepoint
splitting of \S\ref{subsec-stable-basepoint-splitting}. 
\end{proof}

\begin{lemma}\label{lemma-delta-maps-commute-with-projections}
For every unbased space $K$ and spectrum $Y$, the diagram
$$\xymatrix{
{\Sip(\point)\sm \Sip K} \ar[d]_{\Sip \pi}
& {\Sip\Oim Y\sm \Sip K} \ar[l]_-{\Sip\pi\sm 1} \ar[r]^-{\Si q\sm 1}
\ar[d]_{\Sip \delta^+_{Y,K}}
& {\Si\Oi Y\sm \Sip K} \ar[d]^{\Si \delta_{Y,K}}
\\
{\Sip(\point)}
& {\Sip\Oim (Y\sm \Sip K)} \ar[l]^-{\Sip\pi} \ar[r]_-{\Si q}
& {\Si\Oi (Y\sm \Sip K)}
}$$
commutes in the homomotpy category of spectra.  (The map $q$ was
defined in \S\ref{subsec-stable-basepoint-splitting}.)  
\end{lemma}
\begin{proof}
The commutativity of the left-hand square follows from the fact that
$\delta^+$ is a natural transformation, and so commutes with the map
induced by the projection $Y\ra \point$.

The right-hand square is equal to $\Si$ applied to the square
$$\xymatrix{
{(\Oim Y\times K)_+\approx (\Oim Y)_+\sm K_+} \ar[d]_{(\delta^+_{Y,K})_+}
\ar[rr]^-{q\sm 1} 
&& {\Oi Y\sm K_+} \ar[d]^{\delta_{Y,K}}
\\
{(\Oim(Y\sm \Sip K))_+}  \ar[rr]_{q}
&& {\Oi(Y\sm \Sip K)}
}$$
so it suffices to show that this square commutes in the homotopy
category of pointed spaces.  In fact, formal properties of adjunction
show that 
$$q \circ (\delta^+_{Y,K})_+ = \delta^+_{Y,K}= \delta_{Y,K}\circ(q\sm
1).$$ 
\end{proof}

\begin{prop}\label{prop-lambda-is-topological}
Let $K$ be a space and $Y$ a spectrum.  The diagram
$$\xymatrix{
{LY \sm \Sip K} \ar[rr]^-{\lambda^+_Y\sm \id} \ar[d]
&&
{L\Sip\Oim Y\sm \Sip K} \ar[d]^{L(\Sip \delta^+_{Y,K})}
\\
{L(Y\sm \Sip K)} \ar[rr]^-{\lambda^+_{Y\sm \Sip K}}
&& { L\Sip\Oim(Y\sm \Sip K)}
}$$
commutes. (We use the conventions for localization described in
\S\ref{subsec-conventions-on-localization}.) 
\end{prop}

In particular, taking $Y=S$ gives $\lambda^+_{\Sip K}=
L(\Sip\delta_{S,K})\circ (\lambda^+_S\sm \id)$.  

\begin{proof}
This square breaks up into two squares:
$$\xymatrix{
{LY\sm \Sip K} \ar[rr]^-{\lambda_Y\sm \id}  \ar[d]
&& {L\Si\Oi Y\sm \Sip K} \ar[rr]^{L\gamma_{\Oi Y}\sm\id}
\ar[d]_{L\Si\delta_{Y,K}} 
&& {L\Sip\Oim Y\sm \Sip K} \ar[d]^{L\Sip\delta^+_{Y,K}}
\\
{L(Y\sm\Sip K)} \ar[rr]_-{\lambda_{Y\sm \Sip K}}
&& {L\Si\Oi(Y\sm \Sip K)} \ar[rr]_{L\gamma_{\Oi(Y\sm \Sip K)}}
&& {L\Sip\Oim(Y\sm \Sip K)}
}$$
The commutativity of the left-hand square is
\eqref{prop-commutative-lambda-with-rho}, while commutativity of the
right-hand square is proved by applying $L$ to the square of
\eqref{prop-delta-commutes-with-gamma}.
\end{proof}
 
\begin{lemma}\label{lemma-kunneth-product-factorization}
The composite
$$\Oim X\times \Oim Y\xra{\delta^+_{X,\Oim Y}} \Oim(X\sm \Sip\Oim Y)
\xra{\Oim(\id_X\sm \epsilon_Y^+)} \Oim(X\sm Y)$$ is the K\"unneth map,
i.e., the map representing the external product map $X^0K\times
Y^0L \ra (X\sm Y)^0(K\times L)$ in generalized cohomology.  In
particular, for $X=Y=S$, the composite
$$\Oim S\times \Oim S \xra{\delta^+_{S,\Oim S}} \Oim\Sip \Oim S \xra{\Oim
\epsilon_S^+} \Oim S$$
is precisely the K\"unneth product for $S$.
\end{lemma}
\begin{proof}
The K\"unneth map $\Oim X\times \Oim Y\ra \Oim(X\sm Y)$ is characterized
as the adjoint to 
$$\epsilon_X^+\sm \epsilon_Y^+\colon \Sip\Oim X\sm \Sip\Oim Y\ra X\sm
Y.$$
The result follows by factoring $\epsilon^+_X\sm \epsilon^+_Y =
(1\sm \epsilon^+_Y)\circ(\epsilon^+_X\sm 1)$ and taking adjoints.
\end{proof}

As a consequence, we have
\begin{lemma}
\label{lemma-kunneth-product-in-naturality-square}
The diagram
$$\xymatrix{
{LS\sm \Sip \Oim S} \ar[rr]^{\lambda^+_S\sm \id} \ar[d]_{\id\sm \epsilon^+_S}
&& {L\Sip\Oim S\sm \Sip\Oim S} \ar[d]^{L\Sip m}
\\
{LS} \ar[rr]_{\lambda^+_S} 
&& {L\Sip\Oim S}
}$$
commutes, where $m$ denotes the K\"unneth product map $\Oim S\times \Oim
S\ra \Oim S$.  (Recall the conventions described in
\S\ref{subsec-conventions-on-localization} for 
localization.) 
\end{lemma}
\begin{proof}
Consider the diagram
$$\xymatrix{
{LS \sm\Sip\Oim S} \ar[r]^-{\lambda^+_S\sm \id} \ar[d]
& {L\Sip\Oim S\sm \Sip\Oim S} \ar[d]^{L(\Sip\delta^+_{S,\Oim S})}
\\
{L(\Sip\Oim S)} \ar[r]^-{\lambda^+_{\Sip\Oim S}} \ar[d]_{L\epsilon^+_S}
& {L\Sip\Oim(\Sip\Oim S)} \ar[d]^{L\Sip\Oim \epsilon_S^+}
\\
{LS} \ar[r]^{\lambda^+_S} 
& {L\Sip\Oim S}
}$$
The top square is just \eqref{prop-lambda-is-topological} specialized
to $K=\Oim S$ and $Y=S$,  
and so commutes.  The bottom square commutes because $\lambda^+\colon
L\ra L\Sip\Oim$ is a
natural transformation.  The composite of the right-hand vertical maps
is $L\Sip m$, by \eqref{lemma-kunneth-product-factorization}, and
therefore the outer rectangle is the desired square.
\end{proof}

\subsection{Proof of \eqref{thm-log-is-power-operation}}

We must verify for $v$ the
identities (La) and (Lb) of the 
definition of logarithmic element.

\begin{proof}[Proof of \eqref{thm-log-comes-from-logarithmic-elt},
property (La).]  We have that $\epsilon^+=\epsilon\circ q$,
  $\lambda^+=L\gamma\circ \lambda$, and $q\circ \gamma=\id$; thus
  $L\epsilon^+_S\circ \lambda^+_S = \id_S$.
\end{proof}

\begin{proof}[Proof of \eqref{thm-log-comes-from-logarithmic-elt},
property (Lb).]
Apply $R$-homology to the commutative square
of \eqref{lemma-kunneth-product-in-naturality-square} to get
$$\xymatrix{
{\ch{R}(\Oim S)} \ar[r]^-{v\times} \ar[d]_{\tau}
& {\ch{R}(\Oim S\times \Oim S)} \ar[d]^{\circ}
\\
{\ch{R}(\pt)} \ar[r]_{v\cdot}
& {\ch{R}(\Oim S)}
}$$
On the bottom and left, the composite is $x\mapsto \tau(x)\cdot v$.  
On the top and right, the composite is $x\mapsto x\circ v$.
\end{proof}

\section{Level structures and the cocharacter map}
\label{sec-level-structures}

\subsection{Universal deformations}

In this section, we fix a prime $p$ and a height $n\geq1$, a perfect
field $k$ of characteristic $p$, and a height $n$-formal group
$\Gamma_0$ over $k$.  We let $E$ denote the Morava $E$-theory
associated to the universal deformation $\Gamma$ of $\Gamma_0$.

\subsection{Level structures}

For a profinite abelian $p$-group $M$, we write 
$$M^*\defeq \hom^\cts(M,\Q_p/\Z_p),$$
where $\Q_p/\Z_p$ is given the discrete topology, and we write
$M[p^r]$ for the subgroup of $p^r$-torsion elements of $M$.

If $F$ is a formal group over a complete local ring $R$ with maximal
ideal $\m_R$, then $F(\m_R)$ denotes the additive group with
underlying set $\m_R$ and group law given by $F$.  For a discrete abelian
group $A$, a \dfn{homomorphism} $f\colon A\ra F$ means a homomorphism
$A\ra F(\m_R)$ of abelian groups.  The set of such homomorphisms is
denoted $\hom(A,F)$.

\begin{prop}\label{prop-cohomology-of-abelian-group}
Let $A$ be a finite abelian group. 
\begin{enumerate}
\item [(a)]
The $\O$-module $E^0BA^*$ is free and finitely generated over $\O$.
There are natural isomorphisms 
$$E^0(BA^*\times X) \approx E^0BA^*\otimes_{\O} E^0X$$
and
$$R\otimes_{\O}E_0BA^* \approx \hom^\cts_{{\O}\text{-mod}}(E^0BA^*,R).$$
\item [(b)]
Let $i\colon \O\ra R$ be a local homomorphism to a complete local
ring $R$, classifying a deformation $F$.  Then there is a natural isomorphism
$$\hom_{\alg{\O}}^\cts(E^0BA^*, R) \approx \hom(A, F).$$
\end{enumerate}
\end{prop}
\begin{proof}
Part (a) is \cite[5.10 and 5.11]{hopkins-kuhn-ravenel}, and part (b) is
\cite[5.12]{hopkins-kuhn-ravenel}.      
\end{proof}

We set $\O(\hom(A,\Gamma)) \defeq E^0BA^*$; it carries the universal
homomorphism from $A$ to a deformation of $\Gamma_0$.

A homomorphism $f\colon A\ra F$ is called a \dfn{level
structure} if on the formal scheme $F$ over $R$ one has the
inequality of divisiors
$$\sum_{a\in A[p]} [f(a)] \leq F[p],$$
where the left-hand side is over the elements of the $p$-torsion
subgroup of $A$, and
the right-hand side denotes the divisor of the $p$-torsion subgroup of
$F$.    In terms of a coordinate $T$ on $F$, this amounts to
the condition that
$$\prod_{a\in A[p]} (T+_\Gamma T(f(a))) \quad \text{divides}\quad
[p]_F(T)\quad \text{in $R\powser{T}$.}$$
Write $\level(A,F)$ for the set of level
structures; note that by definition it is a subset of
$\hom(A,F)$. 

\begin{prop}\label{prop-cohomology-abelian-groups}
Fix a deformation $F$ of $\Gamma_0$, classified by the homomorphism
$i\colon \O\ra R$.
\begin{enumerate}
\item [(a)]
For each finite abelian group $A$, there exists a complete local ring
$\O(\level(A,\Gamma))$ over $\O$ and natural bijections
$$\hom_{\O-alg}(\O(\level(A,\Gamma)), R) \approx \level(A,F).$$  
The ring $\O(\level(A,\Gamma))$ is a quotient of $\O(\hom(A,\Gamma))$.  
\item [(b)]
If $f\colon A\ra B$ is an inclusion of finite abelian groups, then there is an
evident natural transformation $\level(B,F)\ra
\level(A,F)$.  The map $\O(\level(A,\Gamma))\ra
\O(\level(B,\Gamma))$ classifying the universal example of this
transformation is finite and flat.
\item [(c)] The invariant subring of the evident action of the ring
  $\Aut((\Z/p^r\Z)^n)$ on  $\O(\level((\Z/p^r\N)^n),\Gamma)$ is
  exactly $\O$.
\end{enumerate}
\end{prop}
\begin{proof}
Part (a) is \cite[Proposition
22]{strickland-finite-subgroups-of-formal-groups} or
\cite[10.14]{ando-hopkins-strickland-h-infinity}.  The existence of
the transformation of part (b) is clear.  That the map classifying it
is finite and flat is \cite[Theorem
34(ii)]{strickland-finite-subgroups-of-formal-groups}, while (c) is 
\cite[Theorem 34(iii)]{strickland-finite-subgroups-of-formal-groups}.
\end{proof}

Level structures enter topology in the statement of the character
theorem of \cite{hopkins-kuhn-ravenel}, though this point of view is
not made explicit there.  The most useful references for level
structures in the context of algebraic topology are
\cite[\S7]{strickland-finite-subgroups-of-formal-groups} and
\cite[\S10]{ando-hopkins-strickland-h-infinity}.

\subsection{The cocharacter map}
\label{subsec-cocharacter-map}

Fix $\Lambda\approx \Z_p^n$, so that $\Lambda^*\approx
(\Q_p/\Z_p)^n$.  
Write $D_r \defeq
\O(\level(\Lambda^*[p^r], \Gamma))$, and let $D\defeq \bigcup D_r$.  
The group $GL(\Lambda)$ acts in a natural way on each $D_r$ on the
left, through the finite quotient $GL(\Lambda/p^r\Lambda)$, in such a
way that $D_r^{GL(\Lambda)}\approx \O$.

Let $G$ denote a finite group (not necessarily abelian).  If $M$ is
any profinite abelian $p$-group, we let $G_M\defeq 
\hom^\cts(M,G)/G$, where $G$ acts by conjugation.  In the special case
$M=\Lambda$ it is called the set of
of \dfn{generalized $p$-conjugacy classes}.  There is an evident right
action of $GL(\Lambda)$ on $G_\Lambda$. 

In \cite{hopkins-kuhn-ravenel}, the authors define a character map,
which is a ring homomorphism
$$\chi_G\colon E^0BG \ra \map(G_\Lambda, D)^{GL(\Lambda)}.$$ 
Their theorem \cite[Thm. C]{hopkins-kuhn-ravenel} states that this
becomes an isomorphism after inverting $p$.
It is more convenient for us to use a dual construction, which we call
a \dfn{cocharacter map}.  The cocharacter map $\omega_G\colon G_\Lambda\ra
D\otimes_{\O}E_0BG$ is defined as follows:
An element $x\in G_\Lambda$ is represented by some homomorphism
$f\colon \Lambda/p^r\Lambda\ra G$ for sufficiently large $r$.  There
is a homomorphism $E^0 B\Lambda/p^r\Lambda\twoheadrightarrow
D_r\subseteq D$ classifying the underlying homomorphism of the
universal $\Lambda/p^r\Lambda$-level structure.  Write $\omega_r\in
D\otimes_{\O} \ch{E}B\Lambda/p^r\Lambda$ for the homology class
corresponding to this 
homomorphism, by \eqref{prop-cohomology-of-abelian-group}.  Then we
define  
$$\omega_G(x) \defeq f_*(\omega_r) \in D\otimes_{\O}\ch{E}BG.$$
One checks that this definition does not depend on the choice of
$r$ or $f$, and that the HKR character map is derived from the
cocharacter map: the evaluation of $\chi_G$ at a given $x\in
G_\Lambda$ is given by Kronecker pairing with $\omega_G(x)$.  

Recall that an element $u\in D\otimes_\O \ch{E}X$ is \dfn{grouplike}
if $\pi_*(u)=1$ where $\pi\colon X\ra *$, and if $\Delta_*(u)=u\times
u$, where ``$\times$'' denotes the external product, and $\Delta\colon
X\ra X\times X$ is the diagonal.

\begin{prop}\label{prop-properties-of-cocharacter}\forcepar
\begin{enumerate}
\item [(a)] 
The image of $\omega_G$ is contained in the grouplike
elements of $D\otimes_\O \ch{E}BG$.  

\item [(b)] Under the evident bijection $(G\times H)_\Lambda \approx
G_\Lambda\times H_\Lambda$, we have $\omega_{G\times
H}(x,y)=\omega_G(x)\times \omega_H(y)$.

\item [(c)] The cocharacter map $\omega_G$ is equivariant with respect
to the actions of $GL(\Lambda)$ on $G_\Lambda$ and $D$.

\item [(d)] If $H$ is a subgroup of $G$, and $T\colon \ch{E}BG\ra
\ch{E}BH$ denotes transfer, we have
$$T(\omega_G(x)) = \sum_{gH\in (G/H)^{x(\Lambda)}} \omega_H(x^g),$$
where $x\colon \Lambda\ra G$ is a fixed representative of the
generalized conjugacy class, $x^g(\lambda)\defeq g^{-1}x(\lambda)g$,
and $(G/H)^{x(\Lambda)}$ is the subset of $G/H$ fixed by the image
$x(\Lambda)\subseteq G$.
\end{enumerate}
\end{prop}
\begin{proof}
For (a), it suffices to see that $\omega_r\in D\otimes
\ch{E}B\Lambda/p^r\Lambda$ is grouplike, which is a consequence of
\eqref{prop-cohomology-of-abelian-group}(a) and the
fact that it is dual to a ring homomorphism $E^0B\Lambda/p^r\Lambda\ra
D$. 

The proofs of (b) and (c) are straightforward.

The equation of part (d) is proved in the same way as the transfer
formula for characters \cite[Theorem D]{hopkins-kuhn-ravenel}.
(The statement of the HKR transfer formula directly implies (d) modulo
torsion, which is enough for our purposes.) 
\end{proof}

\subsection{A congruence formula}

Let $i\colon \O\ra \O(\level(A,\Gamma))$ be the standard inclusion,
and let $f\colon A\ra i^*\Gamma$ 
be the tautological level structure.  Let $\I(\level(A,\Gamma))\subset
\O(\level(A,\Gamma))$ denote the ideal generated by the elements
$\{T(f(a))\}_{a\in A}$, where $T$ is any coordinate for $\Gamma$; the ideal
does not depend on the choice of $T$. 

\begin{prop}\label{prop-relation-between-vanishing-ideals}
Suppose $A$ is a finite abelian $p$-group of rank $1\leq r\leq n$.
Then there is
an isomorphism 
$$\O(\level(A,\Gamma))/\I(\level(A,\Gamma)) \approx
\O/(p,\dots,u_{r-1}),$$
and so $\I(\level(A,\Gamma))\cap \O \approx (p,\dots,u_{r-1})$. 
\end{prop}
\begin{proof}
By construction, the $\O$-algebra
$\O(\level(A,\Gamma))/\I(\level(A,\Gamma))$ is 
universal  
for level structures which are trivial homomorphisms.  A deformation
$F$ of $\Gamma_0$ admits at most one such level
structure, and one exists if 
and only if $T^{p^r}$ divides $[p]_{F}(T)$.  Thus $F$ admits
such a trivial level structure if and only if $i(u_k)=0$ for
$k=0,\dots,r-1$, and we conclude that
$\O(\level(A,\Gamma))/\I(\level(A,\Gamma))\approx \O/(p,\dots,u_{r-1})$. 
\end{proof}
\begin{rem}
As a special case of \eqref{prop-relation-between-vanishing-ideals},
we see that
any deformation of $\Gamma_0$ to a ring with $p=0$ 
admits a \emph{unique} $\Z/p$-level structure.  I am indebted to Mike
Hopkins for pointing out this fact to me, which led to the proof of
the congruence \eqref{prop-logarithmic-element-congruence-mod-p}.
\end{rem}

Consider a finite subgroup $V \subseteq \Lambda^*$.  Let
$S\subseteq \Lambda$ denote the kernel of the projection dual to this
inclusion: 
$$0\ra S \ra \Lambda\ra V^*\ra 0.$$
Thus $S$ is an open subgroup of $\Lambda$.  The inclusion $V\subseteq
\Lambda^*$ determines a homomorphism $\O(\level(V,\Gamma))\ra D$,
classifying the restriction of level structures.
\begin{prop}\label{prop-cocharacter-congruence-formula}
Let $x,y\in G_\Lambda$, such that $x|_S=y|_S$ in $G_S$.  Then
$$\omega(x)\equiv \omega(y) \mod \I(\level(V,\Gamma))\cdot
D\otimes_\O \ch{E}BG.$$
\end{prop}
\begin{proof}
Let $\I\defeq \I(\level(V,\Gamma))\cdot D$.
Let $r$ be chosen sufficiently large such that $p^r\Lambda\subseteq
S$, and such that $x$ and 
$y$ are represented by maps $f,g\colon \Lambda/p^r\Lambda\ra G$.  It
suffices to show, in
$$\xymatrix{
{D\otimes_{\O}\ch{E}B\Lambda/p^r\Lambda} \ar[r]^{\pi} 
& {D/\I \otimes_{\O}\ch{E}B\Lambda/p^r\Lambda}  \ar@<-1ex>[d]_{f_*} \ar@<1ex>[d]^{g_*}
& {D/\I \otimes_{\O} \ch{E} BS/p^r\Lambda,} \ar[l]_{j}
\ar[dl]
\\
& {D/\I \otimes\ch{E} BG}
}$$
that $\pi(\omega_r)$ is in the image of $j$; given this, the result
follows, because $x|_S=y|_S$ implies that
$f|_{S/p^r\Lambda}$ and $g|_{S/p^r\Lambda}$ are conjugate by an
element of $G$, and so induce identical maps $\ch{E}BS/p^r\Lambda \ra
\ch{E}BG$.  Dualizing, 
we are 
asking that a dotted arrow exist making the following a commutative
square of rings:
$$\xymatrix{
{E^0B\Lambda/p^r\Lambda} \ar[r]^-{\omega_r} \ar[d] 
& {D} \ar[d]
\\
{E^0BS/p^r\Lambda} \ar@{-->}[r]
& {D/\I}
}$$
The existence of such a dotted arrow is a tautology;  it amounts to a
factorization 
$$\xymatrix{
0 \ar[r]
& {V} \ar[r]
& {\Lambda^*} \ar[r] \ar[d]
& {S^*} \ar[r] \ar@{-->}[dl]
& 0
\\
&& {\Gamma(\m_{D/\I})}
}$$
of a diagram of abelian groups.
\end{proof}

\section{Burnside ring elements}

\subsection{The Morava $E$-theory of $B\Sigma$ and $\Oi S$}

\begin{prop}\label{prop-cohomology-sigma-k-free}
For all $k$, $\ch{E}B\Sigma_k$ is a finitely generated, free $\O$-module.
\end{prop}
\begin{proof}
This is \cite[Theorem 3.3]{strickland-morava-e-theory-of-symmetric}.
\end{proof}

\begin{prop}\label{prop-homology-of-loops-s}
$\ch{E}\Oi S$ is the completion of an infinitely generated free
$\O$-module.  It is flat over $\O$, and thus in particular
is $p$-torsion free.
The union of the images of the maps $\ch{E}\{k\}\times B\Sigma_\ell
\ra \ch{E}\Oi S$ are dense.
\end{prop}
\begin{proof}
Let $f\colon B\Sigma\ra B\Sigma$ denote the map given by $f_k\colon
B\Sigma_{k-1}\ra 
B\Sigma_{k}$.  It is well known that $\Oi S$ is stably equivalent to
$\hocolim(B\Sigma\xra{f}B\Sigma\xra{f}\cdots)$.  
Furthermore, $f_k$ admits a stable retraction $\Sip
B\Sigma_k\ra \Sip B\Sigma_{k-1}$ (the stable homotopy version of a
theorem of Dold \cite{dold-decomposition-theorems}).
From this and
\eqref{prop-cohomology-sigma-k-free}, it is clear that $\pi_k (E\sm
f^{-1}B\Sigma)$ is a free $\O$-module for even $k$, and $0$ for odd
$k$.  Thus $\ch{E}\Oi S$ is the $\m$-adic completion of this free
module.  The flatness result follows from
\eqref{lemma-completion-of-free-module} below.
\end{proof}

\begin{lemma}\label{lemma-completion-of-free-module}
Let $A$ be a Noetherian commutative ring, $I\subseteq A$ an ideal.
Then the $I$-adic completion of any free $A$-module is flat over $A$.
\end{lemma}
\begin{proof}
For a free module on one generator, this is well-known; the usual
proof (e.g., \cite[Prop.~10.14]{atiyah-macdonald-commutative-algebra})
using the 
Artin-Rees lemma generalizes to give the lemma, as we show
below. 

Let $S$ be a set, and define a functor on $A$-modules by $F(M)\defeq
\left(\bigoplus_{s\in S} M\right)^{\sm}_I$.  We claim that
\begin{enumerate}
\item [(i)] $F$ is exact on the full subcategory of finitely generated
  $A$-modules, and
\item [(ii)] the evident map $F(A)\otimes_A M\ra F(M)$ is an
  isomorphism when $M$ is finitely generated.
\end{enumerate}
Then $F(A)\otimes_A{-}$ is exact on the full subcategory of finitely
generated modules, and thus $F(A)$ is flat.

Recall the Artin-Rees lemma
\cite[Thm.~10.11]{atiyah-macdonald-commutative-algebra}: given a
finitely 
generated module $M$ and a submodule $M'$, there exists $c\geq0$ such
that for all $k\geq0$, $I^{k+c}M'\subseteq I^{k}M\cap M' \subseteq
I^{k-c}M'$.  Therefore the same is true when $M$ and $M'$ 
are replaced by $\bigoplus_s M$ and $\bigoplus_s M'$.  This implies
(i), by \cite[Cor.~10.3]{atiyah-macdonald-commutative-algebra}.  To
prove (ii), 
note that the map is an isomorphism if $M$ is free and finitely
generated, and therefore surjective for all finitely generated $M$,
using the exactness result (i).  Since $A$ is Noetherian, (ii)
follows by a $5$-lemma argument.
\end{proof}

\subsection{The cocharacter map for $B\Sigma$ and $\Oi S$}

As in \S\ref{sec-level-structures}, define the profinite abelian group
$\Lambda\defeq \Z_p^n$.
Let $A^+_k(\Lambda)$ denote the set of isomorphism classes of discrete
continuous $\Lambda$-sets which have exactly $k$ elements.  We
identify $A^+_k(\Lambda)\approx  
(\Sigma_k)_\Lambda$: to a generalized conjugacy class $x\colon \Lambda\ra
\Sigma_k$ associate $X=\Sigma_k/\Sigma_{k-1}$, regarded as a
$\Lambda$-set via $x$.

Let $A^+(\Lambda) \defeq
\coprod_k A^+_k(\Lambda)$.  The set $A^+(\Lambda)$ admits the structure
of a semi-ring, with addition and multiplication corresponding to
coproduct and product of sets.  Let $A(\Lambda)$ be the ring obtained
from $A^+(\Lambda)$ by adjoining additive inverses; it is the Burnside
ring of $\Lambda$, isomorphic to the direct limit $A(\Lambda/S)$ where
$S$ ranges over open subgroups of $\Lambda$.

For the following proposition we need the notation introduced in
\S\ref{subsec-structure-maps}.  

\begin{prop}\label{prop-cocharacter-for-symm-groups-etc}
The cocharacter maps \eqref{subsec-cocharacter-map} for symmetric
groups fit together to give a map 
$$\omega^+ \colon A^+(\Lambda) \ra D\otimes_{\O}\ch{E}B\Sigma.$$
It is a homomorphism into the semi-ring of grouplike elements.
That is, for $x,y\in A^+(\Lambda)$,
$$\omega^+(0)=1, \quad\omega^+(1)=[1], \quad
\omega^+(x+y)=\omega^+(x)\cdot \omega^+(y), \quad 
\omega^+(xy)=\omega^+(x)\circ \omega^+(y).$$
If $X$ is a \emph{transitive} $\Lambda$-set, then $\omega^+([X])$ is
also primitive. 

Furthermore, this map extends uniquely to a map
$$\omega\colon A(\Lambda)\ra D\otimes_{\O}\ch{E}\Oi S$$
which is a homomorphism into the ring of grouplike elements in
$D\otimes_{\O} \ch{E}\Oi S$.

We have that for $x\in A(\Lambda)$, 
$$\tau(\omega(x)) = d(x),$$
where $d\colon A(\Lambda)\ra \Z$ is the ring homomorphism
defined by $d([X])=\#(X^\Lambda)$.

The maps $\omega^+$ and $\omega$ are $GL(\Lambda)$-equivariant, and so
$\omega$ induces a map $A(\Lambda)^{GL(\Lambda)}\ra \ch{E}\Oi S$.
\end{prop}
\begin{proof}
That $\omega^+$ lands in the grouplike elements follows from
\eqref{prop-properties-of-cocharacter} (a).   That it is a
homomomorphism of semi-rings follows from the fact that the operations
of sum and product on $A^+(\Lambda)$ are derived, using
\eqref{prop-properties-of-cocharacter} (b), from the maps
$$(\Sigma_k)_\Lambda\times (\Sigma_\ell)_\Lambda \approx
(\Sigma_k\times \Sigma_\ell)_\Lambda \xra{\amalg}
(\Sigma_{k+\ell})_\Lambda$$ 
and 
$$(\Sigma_k)_\Lambda\times (\Sigma_\ell)_\Lambda \approx
(\Sigma_k\times \Sigma_\ell)_\Lambda
\xra{\times}(\Sigma_{k\ell})_\Lambda,$$ 
which are also the origin of the product maps ``$\cdot$'' and
``$\circ$'', as defined in \S\ref{subsec-structure-maps}.   Similarly,
the additive and multiplicative units of $A^+(\Lambda)$ arise as the
unique elements of $(\Sigma_0)_\Lambda$ and $(\Sigma_1)_\Lambda$,
respectively.  

The primitivity of $\omega^+([X])$ when $X$ is transitive follows from
\cite[4.3]{strickland-turner-rational-morava-e-theory}. 

The map $\tau$ is induced by the stable map
$B\Sigma_k\xra{\text{transfer}}B\Sigma_{k-1}\xra{\text{proj}}
\point$.  The transfer formula
\eqref{prop-properties-of-cocharacter}(d) gives
$$\text{transfer}(\omega_{\Sigma_k}(x)) = \sum_{x\in
(\Sigma_k/\Sigma_{k-1})^{x(\Lambda)}} \omega_{\Sigma_{k-1}}(x^g).$$
The element $\omega_G(y)$ is always a grouplike element
\eqref{prop-properties-of-cocharacter}(a), and so goes 
to $1$ under the projection $BG\ra \point$.
Under this projection, the element on the right-hand side of the
equation becomes an integer, 
equal to the size of
$(\Sigma_k/\Sigma_{k-1})^{x(\Lambda)}=X^\Lambda$. 

To extend $\omega^+$ to $\omega$, note that the grouplike elements
of $D\otimes \ch{E}\Oi S$ are invertible in the ``$\cdot$'' product,
so we may set
$$\omega([X]-[Y])\defeq \omega([X])\cdot \omega([Y])^{-1}.$$

The equivariance property follows from
\eqref{prop-properties-of-cocharacter}(c).  
\end{proof}

\subsection{Power operations}\label{subsec-ando-power-operations}

Let $A\subseteq \Lambda^*$ be a finite subgroup of order $p^r$.
Dualizing gives a 
surjective homomorphism $f\colon \Lambda\ra A^*$.  Using this map,  we can
regard $A^*$ as a set with a transitive $\Lambda$-action, and hence an
element in the Burnside ring, denoted $s(A)\in A^+_{p^r}(\Lambda)$.
We define $\psi_A\colon 
E^0(X)\ra D\otimes_\O E^0(X)$ by $\psi_A=\op_{\omega^+(s(A))}$.  
According to the remarks of the previous section, $s(A)$ is both
grouplike and primitive, and thus $\psi_A$ is a ring
homomorphism (though \emph{not} an $\O$-algebra homomorphism).
These operations coincide with the ones constructed by Ando
\cite{ando-power-operations}, though the construction is not
identical, since Ando did not have available to him the fact that the
Morava $E$-theories are commutative $S$-algebras.  Some dicussion of
these operations is given in
\cite{ando-hopkins-strickland-h-infinity}.

\section{Construction of the logarithmic element}
\label{sec-construction-logarithmic-elt}

\subsection{The element defined}\label{subsec-element-defined}

We define a certain element $e\in A(\Lambda)^{GL(\Lambda)}$ as
follows:
$$e= p\sum_{j=0}^n (-1)^j p^{j(j-1)/2} e_j,$$
where 
$$e_j = \frac{1}{p^j}\sum_{\substack{p\Lambda\subseteq S\subseteq
    \Lambda \\ \Lambda/S \approx (\Z/p)^j}} [\Lambda/S].$$
This element $e$ really lives in $A(\Lambda)$ and not just
$A(\Lambda)\otimes \Q$, since $1+j(j-1)/2-j = (j-1)(j-2)/2 \geq0$ when
$j\geq0$.  

In this section, we will prove the following.
\begin{prop}\label{prop-logarithmic-element-congruence-mod-p}
The element $\omega(e)$ is congruent to $1$ modulo $p$ in $\ch{E}\Oi S$.
\end{prop}
\begin{prop}\label{prop-construction-of-logarithmic-elt}
Let $m\in \ch{E}\Oi S$ such that
$1+p\cdot m=\omega(e)$. 
The resulting element 
$$v\defeq \sum_{k=1}^\infty (-1)^{k-1}\frac{p^{k-1}}{k}m^k =
\frac{1}{p}\log\omega(e)$$ 
is the logarithmic element for $E$.
\end{prop}

\subsection{Proof of  the main theorem}
\label{subsec-proof-of-main-theorem}

We can now complete the proof of
\eqref{thm-log-formula-lubin-tate-case}.  By 
\eqref{thm-log-is-power-operation} and
\eqref{thm-log-comes-from-logarithmic-elt}, we have that 
$\ell(\alpha)=\oper{v}(\alpha)$, where $v$ is the logarithmic element for
$E$.   From \S\ref{subsec-structure-maps}(a) and (c) we
have that 
$$\oper{u+u'}(\alpha)=\oper{u}(\alpha)+\oper{u'}(\alpha)\quad \text{and} \quad
\oper{uu'}(\alpha) = \oper{u}(\alpha)\oper{u'}(\alpha),$$
and so 
$$\ell(\alpha)= \sum_{k=1}^\infty (-1)^{k-1}\frac{p^{k-1}}{k}
\oper{m}(\alpha)^k.$$  
The operation $M$ of \eqref{thm-log-formula-lubin-tate-case} is simply
$\op_m$. 

By construction, the operation $M$ satisfies the formula for $1+pM$
given in the statement of \eqref{thm-log-formula-lubin-tate-case}.  We
claim it is the 
unique such operation.  It is clear that the formula characterizes $M$
up to $p$-torsion.  Any operation $E^0\ra D\otimes_\O E^0$
corresponds to an element of $D\otimes_{\O}E^0\Oi E$.  By, e.g.,
\cite[Thm. 1.4]{bendersky-hunton-coalgebraic-ring}, $\ch{E}\Oi E$ is a
free $E_*$-module in even degress, whence $D\otimes_{\O}E^0\Oi E$
is torsion free, 
and thus $M$ must be the unique operation with this property.

\subsection{Congruence for $\omega(e)$}

We use the notation of \S\ref{subsec-ando-power-operations}.  In these
terms, using 
\eqref{prop-cocharacter-for-symm-groups-etc}, we have
$$\omega(e)=\prod_{j=0}^n \biggl(\prod_{\substack{ U\subset
\Lambda^*[p] \\ 
\len{U}=p^j}} \omega(s(U)) \biggr)^{(-1)^j p^{(j-1)(j-2)/2}}.$$

Recall that $\Lambda^*[p]= (\Z/p)^n\subset
(\Q/\Z_{(p)})^n=\Lambda^*$.

\begin{prop}\label{prop-logarithmic-element-congruence-general-case}
In $D\otimes_E \ch{E}\Oi S$ we have the congruence
$$\omega(e) \equiv 1
\mod \I(\level(V,\Gamma))\cdot D\otimes_\O \ch{E}\Oi S,$$
where $V\subset \Lambda^*[p]$ is any subgroup which is isomorphic
to $\Z/p$.  
\end{prop}

\begin{proof}[Reduction of
    \eqref{prop-logarithmic-element-congruence-mod-p} to
    \eqref{prop-logarithmic-element-congruence-general-case}]  
Let $\I=\I(\level(V,\Gamma))$.  
By \eqref{prop-relation-between-vanishing-ideals} we have an inclusion
of short exact sequences 
$$\xymatrix{
{0} \ar[r]
& {p\O} \ar[r] \ar@{>->}[d]
& {\O} \ar[r] \ar@{>->}[d]
& {\O/p\O} \ar[r] \ar@{>->}[d]
& {0}
\\
{0} \ar[r]
& {\I D} \ar[r]
& {D} \ar[r]
& {D/\I D} \ar[r]
& {0.}
}$$
Tensoring with the flat module $\ch{E}\Oi S$
\eqref{prop-homology-of-loops-s} preserves exact sequences and
monomorphisms.  The element $1-\omega(e)\in D\otimes_\O \ch{E}\Oi S$
lives 
in $\ch{E}\Oi S$ by \eqref{prop-cocharacter-for-symm-groups-etc},
and lives in $\I D\otimes_\O \ch{E}\Oi 
S$ by \eqref{prop-logarithmic-element-congruence-general-case}, and so
must be an element of $p\ch{E}\Oi S$.
\end{proof}

\begin{lemma}\label{lemma-congruence-for-cocharacter}
Consider a decomposition $V\oplus V^\perp \approx \Lambda^*[p]$ where
$V\approx \Z/p$;
let $T=\Ker(\Lambda\ra V^*)$.
Given a subgroup $U\subseteq \Lambda^*[p]$, let
$\overline{U}\subseteq V^\perp\subseteq\Lambda^*[p]$ denote the image
of the projection of $U$ to $V^\perp$.  Then in $A^+(T)$,
$$s(U)|_{T} = \begin{cases}
s(\overline{U})|_{T} & \text{if $V\not\subseteq U$,}
\\
p\cdot s(\overline{U})|_{T} & \text{if $V\subseteq U$.}
\end{cases}$$
As a consequence, we obtain the congruences
$$\omega(s(U)) \equiv \begin{cases}
\omega(s(\overline{U})) &\text{if $V\not\subseteq U$,}
\\
\omega(p\cdot s(\overline{U})) & \text{if $V\subseteq U$,}
\end{cases}$$
modulo the ideal $\I(\level(V,\Gamma))\cdot D\otimes_\O \ch{E}\Oi S$. 
\end{lemma}
\begin{proof}
Let $i\colon U\ra \Lambda^*$ denote the given inclusion, and
$j\colon U\ra\Lambda^*$ denote the map factoring through the
projection to $\overline{U}$.  By definiton, $i \equiv j \mod
V$, and so both $i$ and $j$ define the same composite $U\ra
\Lambda^*\ra T^*$.  Dualizing, we see that $i^*|_T=j^*|_T$, viewed as
maps $T\ra\Lambda\ra U^*$.  

If $V\not\subseteq U$, then $j^*|_T$ is surjective, whence
$s(U)|_T=s(\overline{U})|_T$.  If $V\subseteq U$,
then $j^*|_T$ has cokernel isomorphic to $\Z/p$, whence
$s(U)|_T=p\cdot s(\overline{U})|_T$.

The congruences follow immediately from
\eqref{prop-cocharacter-congruence-formula}.  
\end{proof}

\begin{proof}[Proof of
\eqref{prop-logarithmic-element-congruence-general-case}] 
Choose any decomposition $V\oplus V^\perp \approx
\Lambda^*[p]$ with $V\approx \Z/p$, as in the lemma.  Let $d(j)=(-1)^j
p^{(j-1)(j-2)/2}$; 
note that $d(j+1)=-d(j)p^{j-1}$. 
We have 
\begin{align*}
\omega(e) = \prod_{j=0}^n \prod_{\substack{ U\subset \Lambda^*[p] \\
\len{U}=p^j}} \omega(s(U))^{d(j)} 
& = \prod_{j=0}^{n-1} \prod_{\substack{ V\not\subset U \\
\len{U}=p^j}} \omega(s(U))^{d(j)} \cdot \prod_{j=1}^{n} 
\prod_{\substack{V\subset U \\ \len{U}=p^j}}
\omega(s(U))^{d(j)} 
\\
\intertext{which by \eqref{lemma-congruence-for-cocharacter} is
  congruent mod $\I$ to}
& \equiv_{\I} \prod_{j=0}^{n-1}
\prod_{\substack{V\not\subset U \\ \len{U}=p^j}}
\omega(s(\overline{U}))^{d(j)} \cdot 
\prod_{j=1}^{n} \prod_{\substack{V\subset U \\
\len{U}=p^{j}}} \omega(p\cdot s(\overline{U}))^{d(j)}
\\
\intertext{which we reindex according to subgroups of $V^\perp$, to get}
& = \prod_{j=0}^{n-1}
\prod_{\substack{W\subset V^\perp \\ \len{W}=p^j}}
\omega(s(W))^{d(j)p^j} \cdot 
\prod_{j=1}^{n} \prod_{\substack{W\subset V^\perp \\
\len{W}=p^{j-1}}} \omega(s(W))^{d(j)p}
\\
&= \prod_{j=0}^{n-1}\prod_{\substack{W\subset V^\perp \\
\len{W}=p^j}} \omega(s(W))^{d(j)p^j+d(j+1)p}.
\end{align*}
Since the exponents are always $0$, the expression
reduces to $1$.
\end{proof}

\subsection{M\"obius functions and the logarithmic element property}

It remains to show that the element $v$ of
\eqref{prop-logarithmic-element-congruence-general-case} is in fact a
logarithmic element.  To do this, we first show that $e/p\in
A(\Lambda)\otimes 
\Q$ is the idempotent associated to the 
augmentation $d\colon A(\Lambda)\ra \Z$ sending $d([X])=\#X^\Lambda$.

We recall the theory of idempotents in a Burnside ring
\cite{gluck-idempotent-burnside}, in the special case when the group
is \emph{finite abelian}.  Thus, let $G$ be a finite abelian group and
$A(G)$ its 
Burnside ring.  The M\"obius function of 
$G$ is the unique function $\mu_G$ defined on pairs $C\subseteq B$ of
subgroups of $G$, characterized by the property that
$$\sum_{A\subseteq C\subseteq B} \mu_G(C,B) = 
\begin{cases}
  1 &\text{if $A=B$,}
\\
  0 &\text{if $A\neq B$,}
\end{cases}$$
where the sum is over all subgroups $C$ contained in $B$ and
containing $A$.  Then the elements
$$e_A \defeq \sum_{B\subseteq A} \frac{\mu_G(B,A)}{\#(G/B)} [G/B] \in
A(G)\otimes \Q$$
as $A$ ranges over the subgroups of $G$ are the primitive idempotents
of $A(G)\otimes \Q$ \cite[p.~65]{gluck-idempotent-burnside}, and furthermore, 
$$ye_A = d_A(y)e_A,$$
where $d_A\colon A(G)\ra \Z$ is given by $d_A([X])=\#(X^A)$.

Set $\mu_G=\mu_G(0,G)$.  Since the value of $\mu_G(C,B)$ really only
depends on the poset of subgroups of $G$ between $C$ and $B$, we see
that $\mu_G(C,B)=\mu_{B/C}$, and that
$\mu_G$ only depends on the isomorphism class of $G$.  
The following lemma calculates $\mu_A$ for all abelian $A$.
\begin{lemma}\label{lemma-moebius-calculation}
Let $A$ be a finite abelian group.  
\begin{enumerate}
\item[(1)] If $A\approx \prod A_p$ where the $A_p$ are $p$-groups for
distinct primes, then
$\mu_A = \prod \mu_{A_p}$.

\item[(2)] If $A$ is elementary $p$-abelian of rank $j\geq0$, then
$\mu_{A} = (-1)^jp^{j(j-1)/2}$.

\item[(3)] If $A$ is a $p$-group but not elementary $p$-abelian, then
$\mu_{A}=0$.  
\end{enumerate}
\end{lemma}
\begin{proof}
Part (1) is straightforward.

To prove part (2), note that it amounts to 
the identity \cite[Lemma~3.23]{shimura-introduction-automorphic-forms} 
$$\sum_{j=0}^n (-1)^j p^{j(j-1)/2} \gbinom{n}{j}_p = \begin{cases}
1 & \text{if $n=0$,}
\\
0 & \text{if $n>0$,}
\end{cases}$$
where $\gbinom{n}{j}_p = \prod_{i=0}^{j-1}
\frac{(p^{n}-p^i)}{(p^{j}-p^i)}$ is the Gaussian binomial coefficient,
which is the number of elementary abelian subgroups of rank $j$ inside
$(\Z/p)^n$.  

We prove part (3) by induction on the size of $A$.  We have for $A$
nontrivial 
$\mu_A = -\sum_{A\supsetneq B}\mu_B$, the sum taken over proper
subgroups of $A$.  A proper subgroup $B\subsetneq A$ is one 
of two types: (a) it is an elementary abelian $p$-group, or (b) it
isn't.  For (a), such $B$ are exactly the subgroups of $A[p]\subsetneq A$,
the subgroup of $p$-torsion elements, and $\sum_{A[p]\supseteq B}
\mu_B=0$, since $A[p]\neq0$.  For (b), we have $\mu_{B}=0$ by
induction.
\end{proof}

For $r\geq1$, the elements $e_r\defeq e_{\Lambda/p^r\Lambda}\in
A(\Lambda/p^r\Lambda)\otimes \Q$ are idempotents, and the homomorphisms
$A(\Lambda/p^r\Lambda)\otimes \Q\ra
A(\Lambda/p^{r+1}\Lambda)\otimes\Q$ carry $e_r$ to $e_{r+1}$, as can
be seen from the explicit formula for these elements together with
\eqref{lemma-moebius-calculation}.  Thus 
the limiting element $e_\infty\in A(\Lambda)\otimes \Q$ of this
sequence is an idempotent in this ring, with
$[X]e_\infty=d([X])e_\infty$, where $d([X])=\#(X^\Lambda)$.   By
\eqref{lemma-moebius-calculation}, we see that the element $e$ defined
in \S\ref{subsec-element-defined} is equal to $p e_\infty$, and thus
we obtain 
\begin{prop}\label{prop-idempotent-properties}
In $A(\Lambda)$ we have
\begin{enumerate}
\item [(a)] $d(e) = p$, and
\item [(b)] for all $y\in  A(\Lambda)$, $ye=d(y)e$.
\end{enumerate}
\end{prop}

\subsection{Proof of the logarithmic element property}

\begin{lemma}\label{lemma-log-formulas}
Let $\alpha\in D\otimes_\O\ch{E}\Oi S$ be an element of the form
$\alpha=1+p\beta$, and let $\log(\alpha)\defeq \sum_{k\geq1}
(-1)^{k-1}\frac{p^k}{k} 
\beta^{k-1}$.  
\begin{enumerate}
\item [(a)] We have that
$$\tau(\log(\alpha))= \frac{\tau(\alpha)}{\pi(\alpha)}.$$
\item [(b)]
If $w\in D\otimes_\O \ch{E}\Oi S$ is a grouplike element,
then
$$w\circ \log(\alpha)=\log(w\circ \alpha).$$  
\end{enumerate}
\end{lemma}
\begin{proof}
First, note that the operations $\tau$, $\pi$, and $\circ$ on
$D\otimes_\O \ch{E}({-})$ are continuous with respect to the maximal
ideal topology, since they are induced by maps of spectra. 

To prove (a), recall that $\tau$ is
a derivation \eqref{subsec-structure-maps}(m) with respect to the
``$\cdot$'' product, so that
$\tau(\beta^k)=k\tau(\beta)\pi(\beta)^{k-1}$.   
Thus 
\begin{align*}
\tau(\log(\alpha)) 
&= \tau\left( \sum_{k\geq1} (-1)^{k-1} \frac{p^k}{k}\beta^k \right) =
\sum_{k\geq 1} (-1)^{k-1} \frac{p^k}{k}\tau(\beta^k) 
\\
&= \sum_{k\geq 1}(-1)^{k-1}p^k\pi(\beta)^{k-1}\tau(\beta) =
\tau(\beta) (1+p\beta)^{-1}.
\end{align*}

To prove (b), recall that if $w$ is grouplike, then
$w\circ({-})$ is a homomorphism of $D$-algebras
\eqref{subsec-structure-maps}(j).  Thus, 
\begin{align*}
w\circ \log(\alpha)
&= w\circ\left(\sum_{k\geq1}(-1)^{k-1}\frac{p^k}{k}\beta^k \right)
\\
&= \sum_{k\geq1}(-1)^{k-1}\frac{p^k}{k}(w\circ \beta)^k =
\log(1+p(w\circ \beta))=\log(w\circ\alpha).
\end{align*}
\end{proof}

\begin{proof}[Proof of \eqref{prop-construction-of-logarithmic-elt}]
We are going to prove that $v=(1/p)\log \omega(e)$ is a logarithmic element.
We have that
$$\tau(\frac{1}{p}\log\omega(e))
=\frac{1}{p}\frac{\tau(\omega(e))}{\pi(\omega(e))} 
= \frac{1}{p} \frac{d(e)}{1}= \frac{p}{p} = 1,$$  
using \eqref{lemma-log-formulas}(a),
\eqref{prop-cocharacter-for-symm-groups-etc}, and
\eqref{prop-idempotent-properties}(a). 
This proves (La) of the logarithmic element property, that $\tau(v)=1$.

Now we need to prove (Lb); that $x\circ v=\tau(x)v$ for all $x\in
\ch{E}\Oi S$.   
The union of the images of $\ch{E}\{k\}\times B\Sigma_\ell\ra
\ch{E}\Oi S$ is dense, with respect to the maximal ideal topology, by
\eqref{prop-homology-of-loops-s}. 
Thus it suffices to prove (Lb) for those $x$ which are in the image of
one of these maps.  It is enough to do this after faithfully flat base
change to $D$. 
Now, $\ch{E}\Oi S$ is $p$-torsion free, and by the HKR theorem,
$p^{-1}D\otimes \ch{E}B\Sigma_m$ is spanned by elements in the image of the
cocharacter map.  Thus, it suffices to check (Lb) when $x=\omega(y)$
for any $y\in A(\Lambda)$.

So let $y\in A(\Lambda)$.  We have that
\begin{align*}
\omega(y)\circ
\frac{1}{p}\log\omega(e)
& =\frac{1}{p}\log(\omega(y)\circ\omega(e)) 
&& \text{by \eqref{lemma-log-formulas}(b),}
\\
& = \frac{1}{p}\log\omega(ye) = \frac{1}{p}\log\omega(d(y)e) 
&& \text{by \eqref{prop-idempotent-properties}(b),}
\\
&= \frac{1}{p}\log\omega(e)^{d(y)} = \frac{d(y)}{p}\log\omega(e)
\\
&= \tau(\omega(y)) \frac{1}{p}\log\omega(e) 
&& \text{by \eqref{prop-cocharacter-for-symm-groups-etc}.}
\end{align*}
Thus $\omega(y)\circ v = \tau(\omega(y))v$, as desired.
\end{proof}

\section{The logarithm for $K(1)$-local ring spectra}
\label{sec-k1-local}

In this section, we  describe the structure
of $\pi_0 L_{K(1)}\Sip\Oi S$  (completely at an odd
prime, and modulo torsion at the prime $2$), outline its relation
to power operations on $K(1)$-local commutative ring spectra, and give
a proof of \eqref{thm-k1-local-case}.

\subsection{The $p$-adic $K$-theory of some spaces}

We recall results on the $p$-completed $K$-homology of $\Oi S$, due to
\cite{hodgkin-k-theory}.  We have that
$$\ch{K}(B\Sigma_m;\Z_p) \approx \hom_\cts(R\Sigma_m,\Z_p),$$
where $R\Sigma_m$ is the complex representation ring, topologized with
respect to the ideal of representation of virtual dimension $0$.
Furthermore
$$\ch{K}(B\Sigma;\Z_p) \approx
(\Z[\Theta_0,\Theta_1,\Theta_{2},\dots])^{\sm}_p.$$ 
The elements $\Theta_{k}$, $k\geq0$, are characterized implicitly by
Witt polynomials
$$W_k= \sum_{i+j=k} p^i \Theta_{i}^{p^{j}},$$
where $W_k \in K_0(B\Sigma_{p^k};\Z_p)$ is the element
corresponding to the continuous homomorphism $R\Sigma_{p^k}\ra
\Z_p$ defined by evaluation of characters on an
element $g\in \Sigma_{p^k}$ which is a cycle of length $p^k$.  

Thus
$$\ch{K}(\Oi S;\Z_p) \approx
(\Z[\Theta_0^{\pm},\Theta_1,\Theta_2,\dots])^{\sm}_p.$$
Since $W_k\equiv W_0^{p^k}\mod p\ch{K}B\Sigma$, the elements $W_k$ become
invertible in $\ch{K}\Oi S$.  

The cocharacter map $\omega^+\colon A(\Z_p)\ra
\ch{K}(B\Sigma;\Z_p)^{\grplike}$ sends $[\Z_p/p^k]$ to $W_k$.  The
operation corresponding to $W_k$ is the Adams operation $\psi^{p^k}$.  

According to \eqref{prop-construction-of-logarithmic-elt}, the logarithmic
element for $K^{\sm}_p$ is 
$$\frac{1}{p}\log \omega(p[*]-[\Z/p]) = \frac{1}{p}\log
\frac{W_0^p}{W_1} = \frac{1}{p}\log
\frac{1}{1+p\Theta_1/\Theta_0^p}=\sum_{k\geq1}(-1)^{k-1}\frac{p^{k-1}}{k}
\frac{\Theta_1^k}{\Theta_0^{pk}}.$$ 

\subsection{The $K(1)$-local homotopy of $B\Sigma$ and $\Oi S$}

Recall that if $X$ is a spectrum, then there are cofibration sequences
$$\Sigma KO\sm X \ra KO\sm X \ra K\sm X$$
and 
$$LX \ra L(K\sm X) \xra{(\psi^\lambda-1)\sm \id}
L(K\sm X)\quad\text{if $p>2$,}$$ 
and
$$LX \ra L(KO\sm X) \xra{(\psi^\lambda-1)\sm \id}
L(KO\sm X)\quad \text{if $p=2$,}$$
where $\lambda\in \Z_p^\times$ (respectively, $\lambda\in
\Z_2^\times/\{\pm1\}$) is a topological generator.
In particular, we have that
$$\pi_0 L_{K(1)}S = 
\begin{cases}
  \Z_p & \text{if $p>2$,} \\
  \Z_2[x]/(x^2,2x) & \text{if $p=2$,}
\end{cases}$$
where $x$ comes from the non-trivial element of $\pi_1 KO\approx
\Z/2$. 

\begin{prop}
Let $f$ denote either of the Hurewicz maps
$\pi_0 L\Sip B\Sigma_k\ra \ch{K}B\Sigma_k$ or $\pi_0 L\Sip \Oim S \ra
\ch{K}\Oi S$.   The map $f$ is an isomorphism if $p$ is odd, while if
$p=2$ it is surjective, and its kernel is the ideal generated by $x\in
\pi_0LS$. 
\end{prop}
\begin{proof}
This follows from the cofibration sequences mentioned above, together
with the fact that the Adams operations $\psi^\lambda$, for
$\lambda\in \Z_p^\times$, act as the identity map on
$K^0(B\Sigma_k;\Z_p)$. 
\end{proof}

\subsection{The proof of \eqref{thm-k1-local-case}}

Let $R$ be any $K(1)$-local commutative $S$-algebra satisfying the
technical condition described in
\S\ref{subsec-logarithm-for-k1-local}.  By the above 
proposition, there is a natural factorization  $i\colon \ch{K}\Oi S
\ra 
\ch{R}\Oi S$ of the map $\pi_0 L\Sip\Oi S\ra \ch{R}\Oi S$.  Define
power operations in $R$-theory by
$$\psi\defeq \op_{i(W_1)} = \op_{i(\Theta_0^p+p\Theta_1)} \quad\text{and}\quad
\theta\defeq \op_{i(\Theta_1)}.$$
There is an identity $\psi(x)=x^p+p\theta(x)$, and $\psi$ is a ring
homomorphism.  

The derivation of \eqref{thm-k1-local-case} is now straightforward,
since the logarithmic element for $R$ must be the image under $i$ of
the logarithmic element for $K$, since both are the image of the
logarithmic element for $LS$.

\subsection{Exponential maps for $K$-theory and $KO$-theory}

The logarithm maps $\gl_1(K^{\sm}_p)\ra K^{\sm}_p$ and
$\gl_1(KO^{\sm}_p)\ra KO^{\sm}_p$ are seen to be
weak equivalences on $3$-connected covers in the first case, and
$1$-connected covers in the second.  In other words, the logarithm
admits inverse ``exponential'' maps
$$e\colon KSU(X;\Z_p) \ra (1+KSU(X;\Z_p))^\times \quad\text{and}\quad
e\colon KSO(X;\Z_p) \ra (1+KSO(X;\Z_p))^\times,$$
where $KSU({-};\Z_p)$ and $KSO({-};\Z_p)$ denote the cohomology
theories defined by these connective covers.  
Let $\theta_k\colon K^0(X;\Z_p)\ra
K^0(X;\Z_p)$ denote the operation corresponding to the element
$\Theta_k\in \ch{K}B\Sigma_{p^k}$ described above, so that we have
$$\psi^{p^k}(\alpha) = \sum_{i+j=k} p^i\theta_i(\alpha)^{p^j}.$$
\begin{prop}
The exponential maps in $K$ and $KO$ theory are both given by the
formula 
$$e(\alpha)=\prod_{i=0}^\infty \exp\left(\sum_{j=0}^\infty
\frac{\theta_i(\alpha)^{p^j}}{p^j}\right),$$
which converges $p$-adically.  Formally this equals
$\exp\left(\sum_{k=0}^\infty \frac{\psi^{p^k}(\alpha)}{p^k}\right)$. 
\end{prop}

\begin{proof}
We give the proof for $K$-theory; at the end, we indicate the
changes needed for $KO$-theory.

Let $f(T)\defeq
\exp(\sum_{j=0}^\infty \frac{T^{p^j}}{p^j}) \in \Z_{(p)}\powser{T}$ be
the Artin-Hasse exponential.  
We will \emph{define} a map $e$ by
$e(\alpha)\defeq \prod_{i=0}^\infty f(\theta_i(\alpha))$.  We will
show below that this expression
converges for $\alpha\in KSU(X;\Z_p)$ when $X$ is a finite complex.
It is then
straightforward to check that $\ell(e(\alpha))=\alpha$ for any
$\alpha\in KSU(X;\Z_p)$ where $X$ 
is a finite complex.  The representing space for $KSU({-};\Z_p)$ is
$BSU^{\sm}_p$; the set $KSU(BSU^{\sm}_p;\Z_p) = \lim KSU(X;\Z_p)$ as
$X$ ranges over finite subcomplexes, and thus we can verify the
identity $\ell\circ e=\id$ on the universal example, which proves that
$e$ is the desired inverse.

Suppose given a CW-model
$\bigcup X_k =X$.  Let $I_k\defeq \ker[K(X;\Z_p)\ra
K(X_{k-1};\Z_p)]$; the $I_k$'s give a filtration of $K(X;\Z_p)$ by
ideals such that $I_kI_{k'}\subseteq I_{k+k'}$.
We have that $KSU(X;\Z_p)=I_4$.  Since the $\theta_i$ are natural operations,
and $\theta_i(0)=0$, they preserve the ideals $I_k$; in particular,
$\theta_j(I_4)\subseteq I_4$.  Since $X$ is finite, $I_k=0$ for $k$
sufficiently large, and so for $\alpha\in I_4$ each expression
$f(\theta_i(\alpha))$ is actually 
a finite sum, contained in $1+I_4$.  

We now show that $\theta_i(\alpha)\to 0$ as $i\to\infty$ in
the $p$-adic topology.
Since $I_4$ has a finite filtration by the $I_k$'s, it suffices to do
this one filtration quotient at a time.  By the Atiyah-Hirzebruch
spectral sequence, this amounts to a calculation on the
reduced $K$-theory of spheres.  Thus, $I_k/I_{k-1}=0$ if $k$ is odd,
and $\theta_i(\alpha)=p^{i(k/2-1)}\alpha$ for $\alpha\in
I_{k/2}/I_{k/2-1}$ if $k$ is even and $k>0$.  In particular, for
$\alpha\in I_k$, $k\geq 4$, we see 
that the sequence $\theta_i(\alpha) \mod I_{k+1}$ approaches $0$
$p$-adically.  

It follows that $f(\theta_i(\alpha))\to1$ as $i\to\infty$, and
therefore the infinite product converges.

The argument for $KO$-theory is almost the same, except for the
calculations on filtration quotients.  Here the additional
observation is that $\theta_i(\alpha)=0$ for $\alpha\in I_k/I_{k-1}$,
for all $k$ such that $k\geq 2$ and $k\equiv 1,2\mod 8$ (but
\emph{not} when $k=1$). 
\end{proof}

\subsection{Exponential operations of Atiyah-Segal}
In
\cite{atiyah-segal-exponential-isomorphism}, the authors construct
explicit exponential maps on $K$-theory and $KO$-theory completed at
some prime $p$.  Their construction starts with the observation that
on any $\lambda$-ring $R$, the operator
$\Lambda_t\colon R\ra R\powser{t}$ given by
$\Lambda_t(x)=\sum_{i\geq0}\lambda^i(x)t^i$ is exponential.  By
setting $t$ to particular values $\alpha\in \Z_p$, one can sometimes obtain
series which converge $p$-adically, on some subsets of suitable
$p$-adic $\lambda$-rings $R$.  In this way, Atiyah and Segal can piece
together exponential operations, and construct an exponential
isomorphism for $KO$-theory (though not for $K$-theory).  Their
construction involves arbitrary choices, and leads to an operation
which is \emph{not} infinite-loop.

To compare our construction with theirs, we note that in a
$\lambda$-ring we can set $S_t(x)=(\Lambda_{-t}(x))^{-1}$; the
operators $s^i$ defined by $S_t(x)=\sum_{i\geq0}s^i(x)t^i$ correspond
to taking symmetric powers of bundles.  Adams operations are related
to  $S_t$ by the equation
$$S_t(x) = \exp\left[ \sum_{m\geq1} \frac{\psi^m(x)}{m}t^m\right].$$
Thus, our exponential operator is a kind of ``$p$-typicalization'' of
the symmetric powers, evaluated at $t=1$.

\section{The action of Hecke operators on Morava $E$-theory}
\label{sec-hecke-operators}

We give a quick and dirty exposition of a fact which does not seem to
be proved in
the literature, but should be well-known; namely, that the Morava
$E$-theory of a space carries an action by an algebra of Hecke
operators.  

\subsection{Hecke operators}

Let $\Delta$ be a  monoid, and $\Gamma\subset \Delta$ a
sub\emph{group}.  Define
$$\hecke
=\hom_{\Z[\Delta]}(\Z[\Delta/\Gamma],\Z[\Delta/\Gamma])$$ 
where $\Z[\Delta]$ denotes the monoid ring of $\Delta$, and
$\Z[\Delta/\Gamma]$ is the left-$\Z[\Delta]$-module spanned by cosets.
If $M$ 
is a left $\Z[\Delta]$-module, then the $\Gamma$-invariants $M^\Gamma$
are naturally a left $\hecke$-module. 

Consider two examples:
\begin{enumerate}
\item [(a)] The algebra $\hecke=\hecke_n$, associated to $\Delta=\End(\Z^n)\cap
GL_n(\Q)$ and $\Gamma=GL_n(\Z)$.
\item [(b)] The algebra $\hecke=\hecke_{n,p}$, associated to
$\Delta=\End(\Z_p^n)\cap GL_n(\Q_p)$ and $\Gamma=GL_n(\Z_p)$.
\end{enumerate}
In either case, $\hecke$ has a basis which is in one-to-one
correspondence with double cosets $\Gamma\backslash\Delta/\Gamma$; a
double coset $\Gamma x\Gamma$ corresponds to the unique endomorphism
$\tilde{T}_{\Gamma x\Gamma }$ of $\Z[\Delta/\Gamma]$ which sends
$1\mapsto \sum [y\Gamma]$, 
where $y$ ranges over representatives of the finite set $\Gamma
x\Gamma/\Gamma$.  In these terms, $\hecke_n$ is the same as
the Hecke ring for $GL_n$ as described for instance in 
\cite[Ch.~3]{shimura-introduction-automorphic-forms}.  One sees also that
$$\hecke_{n,p}\approx \Z[\tilde{T}_{1,p},\dots,\tilde{T}_{n,p}] \quad
\text{and} \quad \hecke_n \approx \bigotimes_p \hecke_{n,p},$$
where $\tilde{T}_{j,p}$ corresponds to the double coset of the
diagonal matrix which has $p$ in $j$ entries and $1$ in the
other $n-j$ entries.

\subsection{Morava $E$-theory is a module over $\hecke_{n,p}$}

We want to show that the algebra $\hecke_{n,p}$ acts on the functor
$X\mapsto E^0(X)$, where $E$ is a Morava $E$-theory of height $n$.
(Warning: this only agrees up to scalar with the action described in
\S\ref{subsec-interpretation-by-hecke}; see
\S\ref{subsec-renormalized-operators} below.) 
Let $\Lambda=\Z_p^n$.  The right cosets $\Delta/\Gamma$ are in one-to-one
correspondence with open subgroups of $\Lambda$, by $x\Gamma\mapsto
x\Lambda\subseteq \Lambda$.   The sum
$$\sum \omega^+([\Lambda/y\Lambda]) \in D\otimes \ch{E}B\Sigma,$$
where $y$ ranges over representatives of $\Gamma x\Gamma/\Gamma$, is
invariant under the action of $GL(\Lambda)$ on $A^+(\Lambda)$, and so
lives in $\ch{E}B\Sigma$, by
\eqref{prop-cocharacter-for-symm-groups-etc}.  We define 
$$\psi_{\Gamma x\Gamma} \colon E^0X\ra E^0X$$
to be the operation associated to this class.  In terms of the
notation used in \S\ref{subsec-ando-power-operations}, we have
$i\circ \psi_{\Gamma x\Gamma} = \sum \psi_A$, where $i\colon E^0X\ra
D\otimes_\O E^0X$ is the evident inclusion ($D$ is faithfully flat
over $\O$), and the sum is over all
finite subgroups $A\subseteq \Lambda^*$ such that $\ker(\Lambda\ra
A^*)= y\Lambda$ for some $y\in \Gamma x\Gamma$.

\begin{prop}\label{prop-hecke-action-on-lubin-tate}
The assignment $\tilde{T}_{\Gamma x\Gamma}\mapsto \psi_{\Gamma
x\Gamma}$ makes $E^0X$ into an $\hecke_{n,p}$-module.
\end{prop}
This is a statement about compositions of the additive cohomology operations
$\psi_{\Gamma x\Gamma}$.  We will prove it by reducing to results
about the composition of certain \emph{ring} operations, proved in
\cite[App.~B]{ando-hopkins-strickland-h-infinity}.  

\begin{lemma}\label{lemma-ring-ops-hecke}
For each $x\in \Delta$ there is a ring homomorphism
$$\psi_x\colon D\otimes_{\O} E^0X \ra D\otimes_{\O} E^0X$$
natural in $X$, with the following properties:
\begin{enumerate}
\item [(a)] 
If $x\in\Gamma=GL(\Lambda)$, then
$\psi_x$ acts on $D\otimes_\O E^0X$ purely through the
$D$-factor, via the the action of $GL(\Lambda)$ on $D$ described in
\S\ref{subsec-cocharacter-map}.   
\item [(b)] 
Under the inclusion $i\colon E^0X\ra D\otimes_\O E^0X$, we have
$\psi_x\circ i=\psi_A$, where $A\subseteq \Lambda^*$ is the
kernel of the adjoint $x^*\colon \Lambda^*\ra \Lambda^*$ to
$x$. 
\item [(c)] 
If $X$ is a finite product of copies of $\CP^\infty$, then $x\mapsto
\psi_x$ gives an action of 
the monoid $\Delta$ on $D\otimes_\O E^0X$.
\end{enumerate}
\end{lemma}

\begin{proof}
Given $x\in \Delta$, consider the following diagram of formal
groups and level structures over $D$:
$$\xymatrix{
{A} \ar[r]
& {\Lambda^*} \ar[r]^{x^*} \ar[d]_\ell
& {\Lambda^*} \ar[d]^{\ell'}
\\
& {i^*\Gamma} \ar[r]_f
& {j^*\Gamma}
}$$
Here $\Gamma$ is the universal deformation formal group over $\O$,
$i\colon \O\ra D$ is the usual inclusion, $A=\ker x^*$, $f$ is an
isogeny with kernel $\ell(A)$, such that modulo the maximal ideal of
$D$, $f$ reduces to a power of frobenius.  Therefore the codomain of
$f$ is a deformation of $\Gamma_0$ to $D$, classified by a map
$j\colon \O\ra D$.  
There is a commutative diagram
$$\xymatrix{
{D} \ar[dr]_{\chi_{\ell'}}
& {\O=E^0(\point)} \ar[r] \ar[d]^{\psi_A} \ar[l]
& {E^0(X)} \ar[d]^{\psi_A}
\\
& {D=D\otimes_\O E^0(\point)} \ar[r]
& {D\otimes_\O E^0X}
}$$
where $\psi_A$ is the operation associated to $A$ as in
\S\ref{subsec-ando-power-operations}.  On the cohomology of
a point, the map $\psi_A\colon \O\ra D$ equals $j$.  The map
labelled $\chi_{\ell'}$ is the map classifying the pair consisting of
the formal group $j^*\Gamma$ over $D$, and the level structure
$\ell'$.   We define
$$\psi_x\colon D\otimes_\O E^0X \ra D\otimes_\O E^0X\quad
\text{by}\quad x\otimes y\mapsto \chi_{\ell'}(x)\psi_A(y).$$
Properties (a) and (b) are immediate.  Property (c) is proved by the
arguments of \cite[App.~B]{ando-hopkins-strickland-h-infinity} when
$X=\point$ or $X=\CP^\infty$; since the $\psi_x$ act as ring
homomorphisms, 
they are compatible with K\"unneth isomorphisms, and so property (c)
holds for finite products of projective spaces.
\end{proof} 

\begin{proof}[Proof of \eqref{prop-hecke-action-on-lubin-tate}]
We first show that the cohomology operations $\psi_{\Gamma
x\Gamma}$  make 
$E^0X$ into a $\hecke_{n,p}$-module when $X$ 
is a finite product of complex projective spaces,
and hence when $X$ is a finite product of $\CP^\infty$'s.  Since $E$
is a Landweber exact theory, the result of 
\cite[Thm.~4.2]{kashiwabara-hopf-rings} applies to show the desired
result. 

Parts (a) and (c) of \eqref{lemma-ring-ops-hecke} show that the action of
$\Delta$ on $D\otimes_\O E^0X$ descends to an action of $\hecke_{n,p}$
on $(D\otimes_\O E^0X)^\Gamma$.  Then (b) shows that the action of
an operator $\psi_{\Gamma x\Gamma}$ on $E^0X$ as defined above
coincides with  
this action of $\hecke_{n,p}$ under the inclusion $i\colon
E^0X \ra (D\otimes_\O E^0X)^\Gamma$.  
\end{proof}

\subsection{Renormalized operators}
\label{subsec-renormalized-operators}

Setting $T_{j,p}\defeq
(1/p^j)\tilde{T}_{j,p}$, we obtain the operators described in
\S\ref{subsec-interpretation-by-hecke}.  Because of the denominators,
this only gives an action of $\hecke_{n,p}$ on $p^{-1}E^0X$.  We
introduce this apparently awkward renormalization because it coincides
with the usual normalization of Hecke operators acting on classical
modular forms.

\newcommand{\noopsort}[1]{} \newcommand{\printfirst}[2]{#1}
  \newcommand{\singleletter}[1]{#1} \newcommand{\switchargs}[2]{#2#1}
\providecommand{\bysame}{\leavevmode\hbox to3em{\hrulefill}\thinspace}
\providecommand{\MR}{\relax\ifhmode\unskip\space\fi MR }
\providecommand{\MRhref}[2]{%
  \href{http://www.ams.org/mathscinet-getitem?mr=#1}{#2}
}
\providecommand{\href}[2]{#2}


\end{document}